\documentclass{article}
\usepackage[utf8]{inputenc}
\usepackage{amsmath,amssymb,amsthm}
\usepackage[all,cmtip]{xy}
\usepackage{graphicx}
\graphicspath{ {./images/} }
\usepackage{biblatex}
\usepackage{xcolor}
\usepackage{enumerate}
\usepackage{tikz}
\usepackage{authblk}

\makeatletter

\newtheorem{theorem}{Theorem}[section]
\newtheorem*{thm4.1}{Theorem 4.1}

\newtheorem{lemma}[theorem]{Lemma}

\theoremstyle{definition}
\newtheorem{definition}[theorem]{Definition}

\theoremstyle{remark}
\newtheorem{remark}[theorem]{Remark}

\newtheorem{example}[theorem]{Example}

\title{The generalized Harer conjecture for the homology triviality}
\author[1]{Wonjun Chang\thanks{wonjun5@gmail.com}}
\author[1]{Byung Chun Kim\thanks{wizardbc@gmail.com}}
\author[1]{Yongjin Song\thanks{yjsong@inha.ac.kr}}

\affil[1]{\small Department of Mathematics\\Inha University\\Incheon\\22212\\Republic of Korea}
\date{\today}

\begin{document}

\newcommand{\D}{\mathbb{D}}
\newcommand{\R}{\mathbb{R}}
\newcommand{\G}{\mathcal{G}}
\renewcommand{\H}{\mathcal{H}}
\newcommand{\conf}{\operatorname{Conf}}
\newcommand{\homeo}{\operatorname{Homeo}}
\newcommand{\M}{\mathcal{M}}
\newcommand{\B}{\operatorname{B}}
\newcommand{\E}{\operatorname{E}}
\newcommand{\U}{\operatorname{U}}
\newcommand{\supp}{\operatorname{supp}}
\newcommand{\fix}{\operatorname{fix}}
\newcommand{\aut}{\operatorname{Aut}}

\maketitle
\begin{abstract}
The classical Harer conjecture is about the stable homology triviality of the obvious embedding $\phi : B_{2g+2} \hookrightarrow \Gamma_{g}$, which was proved by Song and Tillmann(\cite{ST08}).
The main part of the proof is to show that $\B\phi^{+} : \B B_{\infty}^{+} \rightarrow \B \Gamma_{\infty}^{+}$ induced from $\phi$ is a double loop space map.
In this paper, we give a proof of the generalized Harer conjecture which is about the homology triviality for an $arbitrary$ embedding $\phi : B_{n} \hookrightarrow \Gamma_{g,k}$.
We first show that it suffices to prove it for a $regular$ embedding in which all atomic surfaces are regarded as identical and each atomic twist is a {\it simple twist} interchanging two identical sub-parts of atomic surfaces.
The main strategy of the proof is to show that the map $\Phi : \mathcal{C} \rightarrow \mathcal{S}$ induced by $\B\phi:\conf_n(D)\rightarrow\mathcal{M}_{g,k}$ preserves the actions of the framed little 2-disks operad.
\end{abstract}

\section{Introduction}
\label{sec:intro}

Let $B_{n}$ be Artin's braid group and $\Gamma_{g,k}$ the mapping class group of an orientable surface $S_{g,k}$ with genus $g$ and $k$ boundary components.
There is an obvious embedding $\phi : B_{2g+2} \hookrightarrow \Gamma_{g}, \beta_{i} \mapsto \alpha_{i}$, where $\beta_{i}$ $(1 \leq i \leq 2g+1)$ are the usual generators of $B_{2g+2}$, and $\alpha_{i}$ are the Dehn twists along chain curves $a_{i}$ with $a_{i} \cap a_{i+1} = \{*\}$.
The original Harer conjecture, raised in 1980's, is that
$$ \phi_{*} : H_{*}(B_{\infty}; \mathbb{Z}/2) \rightarrow H_{*}(\Gamma_{\infty}; \mathbb{Z}/2)$$
is trivial.
It has been well-known that it suffice to show that $\phi_{*}$ preserves Araki-Kudo-Dyer-Lashof operations which come from the double loop space structures (\cite{Cohen87}).
We need to show that
$$\B\phi^{+} : \B B_{\infty}^{+} \rightarrow \B \Gamma_{\infty}^{+}$$
is a map of double loop spaces (call this Proposition X).
Proposition X actually implies a very strong fact that $\B\phi^{+}$ is homotopically trivial, because
$\B B_{\infty}^{+} \simeq \Omega^{2}S^{3} \simeq \Omega^{2}\Sigma^{2}S^{1}$ and every double loop space map from $\B B_{\infty}^{+}$ is determined by its restriction $S^{1} \subset \Omega^{2}\Sigma^{2}S^{1}$ up to homotopy.
Note that $\B\Gamma_{\infty}^{+}$ is simply connected by Powell theorem.
Hence Proposition X implies a stronger version of the Harer conjecture :
$$\phi_*:H_i(B_\infty;R)\rightarrow H_i(\Gamma_\infty;R)$$
is trivial for all $i\geq 1$ and any constant coefficient $R$.
This stronger version of Harer conjecture is called the ``homology triviality'' in this paper.

The original Harer conjecture (and its homology triviality) was proved by Song and Tillmann (\cite{ST08}).
They constructed a couple of monoidal 2-categories which correspond to braid groups and mapping class groups, and lifted the embedding $\phi$ to a monoidal 2-functor.

An alternative proof of Proposition X was suggested by Segal and Tillmann (\cite{ST07}).
Their idea is to lift $\phi$ to the map $\Phi : \conf_{2g+2}(D) \rightarrow \mathcal{M}_{g,2}$, where $\conf_{2g+2}(D)$ denotes the configuration space of unordered $2g+2$ distinct points on a disk $D$ and $\mathcal{M}_{g,2}$ denotes the moduli space on surface $S_{g,2}$.
Note that $\conf_{n}(D) \simeq \B B_{n}$ and $\mathcal{M}_{g,2} \simeq \B\Gamma_{g,2}$.
They showed that $\Phi$ is compatible with the natural actions of the framed little 2-disks operad on configuration spaces and moduli spaces, which implies Proposition X.
They interpreted the map $\phi$ as one induced by 2-fold branched covering over a disk with $2g+2$ punctures.
That is, a full Dehn twist on a surface is regarded as a lift, via 2-fold covering, of a half Dehn twist interchanging two points on a disk.

On the other hand, Kim-Song (\cite{KS18}) found the construction of a non-geometric embedding $\phi : B_{n} \hookrightarrow \Gamma_{g,k}$ induced by 3-fold covering, and proved the homology triviality.
They (\cite{KS18},\cite{CKS20}) and Ghaswala-McLeay(\cite{GM20}) independently constructed an infinite family of non-geometric embeddings $\phi_{d}$ induced by $d$-fold ($d \geq 3$) branched coverings.
In section~\ref{sec:d-fold}, we give a proof that every embedding induced by a covering is homologically trivial (Theorem~\ref{thm:harer}).
We may have even more general result.
The following theorem is the main result of this paper which is most generalized form of Harer conjecture.

\begin{thm4.1}
For every embedding $\phi : B_{n} \hookrightarrow \Gamma_{g,k}$,
$$\phi_{*}:H_{i}(B_{\infty}; R) \rightarrow H_{i}(\Gamma_{\infty}; R)$$
is trivial for all $i\geq 1$ and any coefficient $R$.
\end{thm4.1}

To every embedding $\phi : B_{n} \hookrightarrow \Gamma_{g,k}$, there is a corresponding space map $\Phi : \conf_{n}(D) \rightarrow \mathcal{M}_{g,k}$ by Earle-Eells theory.
We show that the map between two algebras induced by $\Phi$ preserves the actions of the framed little 2-disks operad.

Let $\phi : B_{n} \hookrightarrow \Gamma_{g,k}, \beta_{i} \mapsto [\tau_{i}]$, where each $\beta_{i}$ is the standard generator interchanging the $i$-th and $(i+1)$-st points.
Each $\tau_{i}$ is taken to be a self-homeomorphism of the surface satisfying the braid relation
$[\tau_{i}][\tau_{i+1}][\tau_{i}] = [\tau_{i+1}][\tau_{i}][\tau_{i+1}]$.
$\tau_{i}$ may be regarded as one isotopic to a self-homeomorphism interchanging two identical parts of the surface $S_{g,b}$ which are points or subsurfaces, not necessarily connected, because a standard generator of braid group is, in principle, characterized as interchanging two points by a half twist.

For a self-homeomorphism $f$ of a surface, the support $\supp(f)$ is defined to be the set of non-fixed points.
For an embedding $\phi : B_{n} \hookrightarrow \Gamma_{g,k}$ we may regard $S_{g,k}$ as the union of all $\overline{\supp([\tau_{i}])} = T_{i}$, for $i = 1,\cdots, n-1$.
Recall that by Harer-Ivanov theorem, the homology of mapping class group $\Gamma_{g,k}$ is, in a stable range, independent of $g$ and $k$.
Each $T_{i}$ is a subsurface of $S_{g,k}$ and $\tau_{i}$ fixes the boundary components of $T_{i}$ pointwise and moves all interior points of $T_{i}$ while interchanging two identical subsurfaces (or points) $I_{i}$ and $I_{i+1}$ of $T_{i}$ by a `simple twist' (Definition~\ref{def:half_twist}).
We may also assume that all $T_i$ are homeomorphic to each other (Lemma~\ref{lem:supp}) and connected.
We define such an embedding to be a regular embedding (Definition~\ref{def:regular}).
For the proof of Theorem~\ref{thm:GeneralHarer}, it suffice to show the homology triviality for regular embeddings.

In order to show that the algebra map induced by $\Phi : \conf_{n}(D) \rightarrow \mathcal{M}_{g,k}$ preserves the actions of the framed little 2-disks operad, we need to figure out the formula for the number of boundary components of the surface, because the action of the operad on the moduli spaces is a kind of `capping' the holes of the surfaces by the framed little 2-disks operad.

In section 2, we recall the construction of the embedding $\phi_{d} : B_{n} \rightarrow \Gamma_{g,k}$ induced by $d$-fold covering over a disk, and show that it satisfies the homology triviality by constructing the algebra map over the framed little 2-disks operad.

In section 3, we introduce the definition of regular embedding and give a few examples of non-geometric embeddings which are all regular embeddings.

In section 4, we prove the main theorem in four steps.
In step I, we reduce the theorem to the case where the embedding is regular.
In step II, we detach $I_{i}$ and $I_{i+1}$ from each $T_{i}$ (note that $I_{i+1}\subset T_{i}\cap T_{i+1}$) and let the holes (the boundary components of $I_{i}$ and $I_{i+1}$) collapse to points.
Because all $I_{i}$ are identical, we may consider only $I_{1}$.
Let $\widetilde T_{1}$ be the new surface obtained by the collapse of the boundary components of $I_{1}$ and $I_{2}$.
Let $\widetilde\tau_{1}$ be the restriction of $\tau_{1}$ onto $\widetilde T_{1}$.
Then $\widetilde\tau_{1}$ acts on the $j$-th component of $\widetilde T_{1}$ as interchanging two points $p_{1,j}$ and $p_{2,j}$ by so-called simple twist.
In this step we first consider only the case where $\widetilde T_{1}$ is connected.
In step III, we show that $\widetilde\tau_{1}$ can be represented by a self-functor of a groupoid and the number of boundary components of $\widetilde T_{1}$ equals 1 or 2.
That is, $\widetilde\tau_{1}$ is equivalent to the lift of half Dehn twist by $d$-fold covering over a disk with two points.
Therefore, by the same reason as in section 2, the homology triviality holds.
In the final step, we complete the proof of the theorem by counting the case where $\widetilde T_1$ is disconnected.

\section{The homology triviality of the embedding $\phi_d$}
\label{sec:d-fold}

Let $B_{n}$ be the braid group with the standard generators $\beta_{i}$ ($1 \leq i \leq n-1$).

\begin{definition}
  An embedding $\phi : B_{n} \hookrightarrow \Gamma_{g,k}$ is said to be $geometric$ if it maps each $\beta_{i}$ to a Dehn twist in $\Gamma_{g,k}$.
\end{definition}

There is an obvious geometric embedding $B_{2g+2} \hookrightarrow \Gamma_{g}$ mapping $\beta_{i}$ to the Humphries generators $\alpha_{i}$ which are mapping classes of consecutive Dehn twists.
In \cite{ST08}, Segal and Tillmann viewed this embedding as being induced by a 2-fold branched covering over a disk.
For $d \geq 3$, the embeddings $\phi_{d}$, induced by $d$-fold coverings were constructed by Kim-Song, Ghaswala-McLeay (\cite{CKS20}, \cite{GM20}, \cite{CS20}).
Each $\phi_{d}(\beta_{i})$ is proved to be a product of $d-1$ Dehn twists, and hence $\phi_{d}$ are all non-geometric for $d \geq 3$.
This gives another answer to the question of Wajnryb \cite{Wajnryb06} about the existence of non-geometric embedding.

Let $p : S_{g,b}^{(n)} \rightarrow S_{0,1}^{(n)}$ be a cyclic $d$-fold branched covering with $n$ branch points $p_{1}, \ldots , p_{n}$.
A cyclic branched covering means that around the branch points, $p$ maps as $z \mapsto z^d$.

We may regard $S_{0,1}^{(n)}$ as a full subcategory (subgroupoid) $D$ of the fundamental groupoid $\Pi_1(S_{0,1})$ whose set of objects is denoted by $\{p_{0}, p_{1}, \ldots , p_{n}, p_{n+1}\}$ where $p_{0}$ and $p_{n+1}$ are on the boundary component.
Then the morphisms of $D$ is generated by $\{e_0,e_1,\ldots,e_{n}\}$ where $e_i$ represents the homotopy class of a path from $p_i$ to $p_{i+1}$.
Similarly, we can regard $S_{g,b}^{(n)}$ as a full subcategory $E$ of $\Pi_1(S_{g,b})$ whose set of objects is $p^{-1}(p_0)\cup p^{-1}(p_{n+1})\cup\{p_1,\ldots,p_n\}$,
and the morphisms of $E$ is generated by the homotopy classes of paths $e_{i,j}$ $(0\leq i\leq n, 1\leq j\leq d)$ corresponding to $e_i$ of $D$ through cyclic $d$-fold covering map $p$.
For more details, see \cite{KS18} and \cite{CKS20}.

\begin{remark}
\label{rmk:ribbon}
The above groupoid argument can also be interpreted in terms of ribbon graphs.
Let $G = (V, E)$ be the graph with the set of vertices $V = \{p_{1}, \cdots , p_{n}\}$ and the of edges $E = \{e_{i,j}\}_{1 \leq i \leq n-1, 1 \leq j \leq d}$ where $e_{i,j}$ (abuse of notation) is an edge corresponding to the homotopy class of the path $e_{i,j}$.
We view each $e_{i,j}$ to be an edge between $p_{i}$ and $p_{i+1}$.
The midpoint of the edge $e_{i,j}$ splits the edge into two parts, which are called the half edges.
A half edge incident to $p_{i}$ is denoted by $e_{i,j}^{-}$, and a half edge incident to $p_{i+1}$ is denoted by $e_{i,j}^{+}$.

Let $\mathcal G$ be the ribbon graph obtained from $G$ with the set of half-edges $H = \{e_{i,j}^-, e_{i,j}^+\}_{1 \leq i \leq n-1, 1 \leq j \leq d}$.
The cyclic ordering on edges at $p_{1}$ is given by $(e_{1,1}^{+}, e_{1,2}^{+}, \cdots , e_{1,d}^{+})$, and at $p_{n}$, it is given by $(e_{n-1,1}^{-}, e_{n-1,2}^{-}, \cdots , e_{n-1, d}^{-})$.
For $i \not\in \{1, n-1\}$, the cyclic ordering at $p_{i}$ is given by
$(e_{i-1,1}^{-}, e_{i,1}^{+}, e_{i-1,2}^{-}, e_{i, 2}^{+}, \cdots , e_{i-1, d}^{-}, e_{i,d}^{+})$.
Note that the surface $S_{g,b}$ can be obtained as a geometric realization of this ribbon graph.
\end{remark}
\ \\

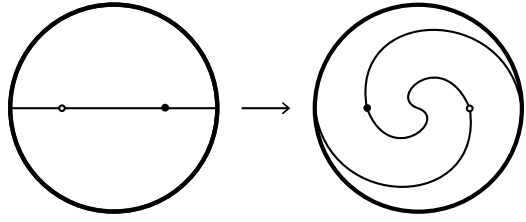
\begin{figure}[h]
    \centering
    \definecolor{cffffff}{RGB}{255,255,255}
    \def \globalscale {0.8}
    \begin{tikzpicture}[y=0.80pt, x=0.80pt, yscale=-\globalscale, xscale=\globalscale, inner sep=0pt, outer sep=0pt]
      \path[draw=black,line width=1.600pt] (100.8000,147.6000) .. controls (100.8000,113.8000) and (128.2000,86.4000) .. (162.0000,86.4000) .. controls (195.8000,86.4000) and (223.2000,113.8000) .. (223.2000,147.6000) .. controls (223.2000,181.4000) and (195.8000,208.8000) .. (162.0000,208.8000) .. controls (128.2000,208.8000) and (100.8000,181.4000) .. (100.8000,147.6000) -- cycle;

      \path[draw=black,line width=0.800pt] (100.8000,147.6000) -- (223.2000,147.6000);

          \path[draw=black,line width=0.800pt] (265.1090,147.6000) .. controls (262.7210,147.6000) and (237.6000,147.6000) .. (237.6000,147.6000);

          \path[fill=black] (262.2590,150.9870) -- (265.9260,147.9870) -- (266.3990,147.6000) -- (265.9260,147.2130) -- (262.2590,144.2130) .. controls (262.0450,144.0380) and (261.7300,144.0700) .. (261.5550,144.2830) .. controls (261.3810,144.4970) and (261.4120,144.8120) .. (261.6260,144.9870) -- (265.2920,147.9870) -- (265.2920,147.2130) -- (261.6260,150.2130) .. controls (261.4120,150.3880) and (261.3810,150.7030) .. (261.5550,150.9170) .. controls (261.7300,151.1300) and (262.0450,151.1620) .. (262.2590,150.9870) -- cycle;

      \path[draw=black,line width=1.600pt] (280.8000,147.6000) .. controls (280.8000,181.4000) and (308.2000,208.8000) .. (342.0000,208.8000) .. controls (375.8000,208.8000) and (403.2000,181.4000) .. (403.2000,147.6000) .. controls (403.2000,113.8000) and (375.8000,86.4000) .. (342.0000,86.4000) .. controls (308.2000,86.4000) and (280.8000,113.8000) .. (280.8000,147.6000) -- cycle;

      \path[draw=black,line width=0.800pt] (280.8000,147.4650) .. controls (288.4100,210.4930) and (381.6550,208.8000) .. (372.3730,147.6000);

      \path[draw=black,line width=0.800pt] (311.6270,147.6410) .. controls (325.4200,183.9530) and (360.8320,154.1170) .. (341.7980,147.6000) .. controls (322.2860,140.9190) and (356.8510,110.3890) .. (372.3730,147.6000);

      \path[draw=black,line width=0.800pt] (403.2000,147.7760) .. controls (395.5910,84.7479) and (300.5680,86.7958) .. (311.6270,147.6410);

      \path[draw=black,fill=cffffff,line width=0.800pt] (370.6900,147.7350) .. controls (370.6900,148.6650) and (371.4440,149.4190) .. (372.3730,149.4190) .. controls (373.3030,149.4190) and (374.0570,148.6650) .. (374.0570,147.7350) .. controls (374.0570,146.8060) and (373.3030,146.0520) .. (372.3730,146.0520) .. controls (371.4440,146.0520) and (370.6900,146.8060) .. (370.6900,147.7350) -- cycle;

      \path[draw=black,fill=black,line width=0.800pt] (310.1690,147.4650) .. controls (310.1690,148.3720) and (310.9040,149.1070) .. (311.8110,149.1070) .. controls (312.7180,149.1070) and (313.4540,148.3720) .. (313.4540,147.4650) .. controls (313.4540,146.5580) and (312.7180,145.8230) .. (311.8110,145.8230) .. controls (310.9040,145.8230) and (310.1690,146.5580) .. (310.1690,147.4650) -- cycle;

      \path[draw=black,line width=1.600pt] (100.8000,147.4240) .. controls (100.8000,113.6240) and (128.2000,86.2236) .. (162.0000,86.2236) .. controls (195.8000,86.2236) and (223.2000,113.6240) .. (223.2000,147.4240) .. controls (223.2000,181.2230) and (195.8000,208.6240) .. (162.0000,208.6240) .. controls (128.2000,208.6240) and (100.8000,181.2230) .. (100.8000,147.4240) -- cycle;

      \path[draw=black,fill=black,line width=0.800pt] (190.6900,147.2880) .. controls (190.6900,146.3590) and (191.4440,145.6050) .. (192.3730,145.6050) .. controls (193.3030,145.6050) and (194.0570,146.3590) .. (194.0570,147.2880) .. controls (194.0570,148.2180) and (193.3030,148.9720) .. (192.3730,148.9720) .. controls (191.4440,148.9720) and (190.6900,148.2180) .. (190.6900,147.2880) -- cycle;

      \path[draw=black,fill=cffffff,line width=0.800pt] (129.7540,147.6000) .. controls (129.7540,146.6930) and (130.4900,145.9580) .. (131.3960,145.9580) .. controls (132.3030,145.9580) and (133.0390,146.6930) .. (133.0390,147.6000) .. controls (133.0390,148.5070) and (132.3030,149.2420) .. (131.3960,149.2420) .. controls (130.4900,149.2420) and (129.7540,148.5070) .. (129.7540,147.6000) -- cycle;

    \end{tikzpicture}

    \caption{Half Dehn twist, $re^{i\theta}\mapsto re^{i(\theta-2\pi r)}$ where $r\in[0,1]$ and $\theta\in[0,2\pi]$.}
    \label{fig:half_Dehn}
\end{figure}
The generators $\beta_{i}$ of $B_n$ are the half Dehn twists around the two points $p_{i}$ and $p_{i+1}$ in $S_{0,1}^{n}$ as in Figure~\ref{fig:half_Dehn}.
The half Dehn twist $\beta_i$ on $S_{0,1}^{(n)}$ can be expressed as a self-functor of the groupoid $D$ as follows :

$$
    \beta_i:\left\{
    \begin{array}{rcl}
      p_i & \mapsto & p_{i+1} \\
      p_{i+1} & \mapsto & p_i \\
      \\
      e_{i-1} & \mapsto & e_{i-1} \cdot e_{i} \\
      e_{i} & \mapsto & e_{i}^{-1} \\
      e_{i+1} & \mapsto & e_{i} \cdot e_{i+1}
    \end{array}\right.
$$
for each $i=1,\ldots, n-1$.
It fixes points and edges that do not appear in the list.

We can see that for $d$ odd, $\beta_{i}$ is lifted to a $\frac{1}{2d}$-Dehn twist, and for $d$ even, $\beta_{i}$ is lifted to a couple of $\frac{1}{d}$-Dehn twists.
Let $\widetilde\beta_i$ be the self-homeomorphism of $S_{g,b}^{(n)}$ satisfying the following commutative diagram:
$$
\xymatrix{
S_{g,b}^{(n)} \ar[d]_p \ar[r]^{\widetilde\beta_i} & S_{g,b}^{(n)}\ar[d]^p\\
S_{0,1}^{(n)} \ar[r]^{\beta_i} & S_{0,1}^{(n)}.
}
$$
As a self-functor of the groupoid $E$, $\widetilde\beta_i$ acts on points and edges as follows (\cite{CKS20}):
$$
\widetilde\beta_i:\left\{
\begin{array}{rcl}
  p_i & \mapsto & p_{i+1} \\
  p_{i+1} & \mapsto & p_i \\
  \\
  e_{i-1, j} & \mapsto & e_{i-1, j} \cdot e_{i, j+1} \\
  e_{i,j} & \mapsto & e_{i, j+1}^{-1} \\
  e_{i+1, j} & \mapsto & e_{i,j} \cdot e_{i+1, j}.
\end{array}\right.
$$

It was shown in \cite{KS18} and \cite{CKS20} that the homomorphism $\phi_{d} : B_{n} \rightarrow \Gamma_{g,b}$, $\beta_{i} \mapsto \widetilde{\beta}_{i}$, induced by $d$-fold cyclic branched covering over a disk with $n$ branch points, is well-defined and injective by the Birman-Hilden theory.

Now, let us show the homology triviality of the embedding $\phi_{d}$.
The proof is similar to the case $d=3$ given in \cite{KS18}.
The main part is to show that the action of the framed little 2-disks operad is preserved by the map induced from $\phi_{d}$.
This will prove Proposition X for $\phi_d$ which implies the homology triviality.

Let $D$ be the unit closed disk $\{z \in \mathbb{C} \mid |z|  \leq 1\}$, and $D_{1}, \cdots ,D_{k}$ be $k$ disjoint, ordered copies of $D$.
Define $\mathcal{D}_{k} = \{f : D_{1} \sqcup \cdots \sqcup D_{k} \rightarrow D \mid \left.f\right|_{D_{i}}(x) = \alpha_{i} x + \beta_{i}, \text{where } \alpha_{i} \in \mathbb{C}^{\times}, \beta_{i} \in \mathbb{C}, f(D_i) \cap f(D_j) = \phi \text{ for } i \neq j \}$, the space of smooth embeddings of $k$ disjoint ordered copies of $D$, that restrict on each $D_{i}$ by the composition of a translation and a multiplication by an element of $\mathbb{C}^{\times}$.

The spaces $\{D_{n}\}_{n \geq 0}$ form an operad $\mathcal{D}$, called the framed little 2-disks operad, with the structure map
$$ \gamma : \mathcal{D}_{k} \times (\mathcal{D}_{m_{1}} \times \cdots \mathcal{D}_{m_{k}}) \rightarrow \mathcal{D}_{\Sigma m_{i}} $$
given by composition of embeddings. (\cite{ST07}, \cite{KS18})

Let $\conf_{n}(D)$ denote the configuration space of unordered $n$-tuples of distinct points in the interior of $D$.
Recall that $\conf_n(D)\simeq\B B_n$.
Put $\mathcal{C}_{m} = \conf_{dm}(D)$ $(m \geq 0)$, and $\mathcal{C} = \{\mathcal{C}_{m}\}_{m \geq 0}$. Then $\mathcal C$ is a $\mathcal D$-algebra with the action
$$ \gamma_{\mathcal{C}} : \mathcal{D}_{k} \times (\mathcal{C}_{m_{1}} \times \cdots \mathcal{C}_{m_{k}}) \rightarrow \mathcal{C}_{\Sigma m_{i}} $$
defined by $(f; a_{1}, a_{2}, \cdots , a_{k}) \mapsto f(a_{1} \cup a_{2} \cup \cdots \cup a_{k})$.

Let $\mathcal{M}_{g,b}$ be the moduli space of Riemann surfaces of genus $g$ with $b$ boundary components.
Recall that $\mathcal M_{g,b}\simeq\B\Gamma_{g,b}$.
For the embedding $\phi_{d} : B_{dm} \hookrightarrow \Gamma_{g,d}$,
we have $g = g(d,m) = (d^{2}m - md - 2d + 2)/2$ by the Riemann-Hurwitz formula.
Let $\mathcal{S} = \{\mathcal{S}_{m}\}_{m \geq 0}$, where

$$
    \mathcal{S}_{m}:=\left\{
    \begin{array}{lcl}
      \mathcal{M}_{g(d,m), d} & \text{ for } & m \geq 1 \\
      \mathcal{M}_{0,1} \sqcup \cdots \sqcup \mathcal{M}_{0,1} (d \text{ times}) & \text{ for } & m = 0.
    \end{array}\right.
$$
$\mathcal S$ is obviously a $\mathcal D$-algebra.
Then the embedding $\phi_d:B_{dm}\hookrightarrow\Gamma_{g(d,m),d}$ may be lifted to the map $\Phi_d:\mathcal C\rightarrow\mathcal S$ through the classifying space functor.

\begin{theorem}
$\Phi_d : \mathcal{C} \rightarrow \mathcal{S}$
is a map of $\mathcal{D}$-algebras.
\end{theorem}

\begin{figure}[h]
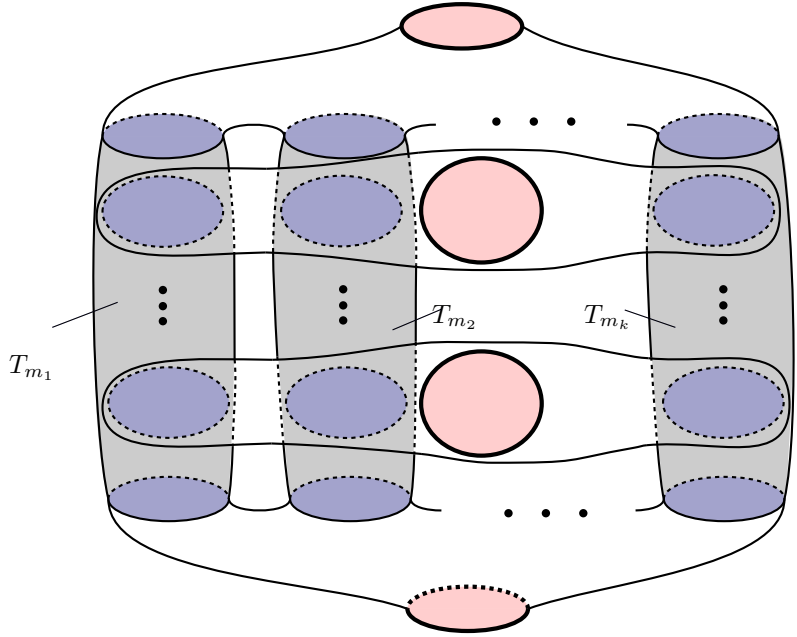

    \centering
    \definecolor{cff0000}{RGB}{255,0,0}
    \definecolor{c0000ff}{RGB}{0,0,255}
    \def \globalscale {0.8}

    \caption{The action of the framed little 2-disks operad on $\mathcal S$.}
    \label{fig:operad}
\end{figure}

\begin{proof}

Each surface $T$ in $\mathcal{S}_{m}$ has $d$ ordered parametrized boundary circles.
For $f \in \mathcal{D}_{k}$, and $j \in \{1,\cdots,d\}$,
let $(P_{f})_{j} = D \backslash f(D \cup \cdots \cup D) = S_{0, k+1}$. 
We think of $(P_{f})_{j}$ as a Riemann surface with $k+1$ parametrized boundary components.
For $f \in \mathcal{D}_{k}$, and $T_{m_{i}} \in  \mathcal{S}_{m_{i}}$, let us place the surfaces $T_{m_{i}}$ from left to right as shown in the Figure~\ref{fig:operad}.
For each $m_i$, we glue the boundary components of $(P_{f})_{i}$ to the boundary circles of $T_{m_i}$.
The resulting surface has extra genus, $(k-1)(d-1)$, that is, the total genus of the resulting surface is $g(d,m_{1}) + \cdots + g(d,m_{k}) + (k-1)(d-1)$ which equals $g(d,\Sigma m_{i})$. 
This gives us the action of the framed little 2-disks operad,
$$ \gamma_{\mathcal{S}} : \mathcal{D}_{k} \times (\mathcal{S}_{m_{1}} \times \cdots \mathcal{S}_{m_{k}}) \rightarrow \mathcal{S}_{\Sigma m_{i}}. $$
Hence, $\Phi_d$ is a $\mathcal{D}$-algebra map.
\end{proof}

Let $X^{+}$ be the Quillen's plus construction, and $\Sigma X$ be its (reduced) suspension.
Since $\phi$ induces a map $\Phi_d : \mathcal{C} \rightarrow \mathcal{S}$ of $\mathcal{D}$-algebras over the little 2-disks operad, $\Phi_d$ induces a map $\Phi_d^{+}$ between their group completions, which is a map of double loop spaces.
We have

$$\Phi_d^{+} : \Omega \B(\amalg_{m > 0} \B B_{dm}) \simeq \mathbb{Z} \times \B B_{\infty} ^{+} \rightarrow \Omega \B(\amalg_{m > 0} \B \Gamma_{g(d,m),d}) \simeq \mathbb{Z} \times \B\Gamma_{\infty} ^+.$$
We know that $\B B_{\infty}^{+} \simeq \Omega^{2}S^{3}$ and $\B\Gamma_{\infty}^{+}$ is an infinite loop space.
The embedding $\phi_{d}$ induces a map on double loop spaces $\B\phi_{d}^{+} : \B B_{\infty}^{+} \rightarrow \B\Gamma_{\infty}^{+}$.
This map is completely determined (upto homotopy) by its restriction to $S^{1} \subset \B B_{\infty}^{+} \simeq \Omega^{2}S^{3} \simeq \Omega^{2} \Sigma^{2} S^{1}$.
By Powell Theorem, the mapping class group is perfect, that is, $\B \Gamma_{\infty}^{+}$ is simply connected.
Hence $\Phi_d^{+}$ is null-homotopic on the connected component, which obviously implies the homology triviality.

\begin{theorem}
  \label{thm:harer}
  Let $\phi_{d} : B_{n} \hookrightarrow \Gamma_{g,b}$ be the embedding induced by $d$-fold branched covering ($d\geq 2$).
  Then for all $i \geq 1$, 
  $$(\phi_{d})_{*} : H_{i}(B_{\infty}; R) \rightarrow H_{i}(\Gamma_{\infty}; R)$$
  is trivial in any constant coefficient $R$.
\end{theorem}

\section{Regular embeddings}
\label{sec:regular}

For a self-homeomorphism $f$ of $S_{g,b}$, the {\it support} $\supp(f)$ is defined to be the set of non-fixed points of $f$.

Let $\phi:B_n\hookrightarrow\Gamma_{g,k}$ be an embedding
and $\beta_1,\ldots,\beta_{n-1}$ be the standard generators of $B_n$.
Since $\phi(\beta_i)$ satisfy the braid relations,
we may take the representative self-homeomorphism $\tau_i$ of $\phi(\beta_i)$
as a kind of {\it simple twist} interchanging two identical subsurfaces (or points) $I_i$ and $I_{i+1}$ on $T_i=\overline{\supp(\tau_i)}$.
Since all $\tau_i$ will be taken to be conjugate to each other,
we may assume that all $T_i$ are identical by the following lemma.

\begin{lemma}
\label{lem:supp}
Let $f$ be a self-map of $S_{g,b}$ and let $h$ be a self-homeomorphism of $S_{g,b}$. Then we have
$$\supp(f)\cong h^{-1}\left(\supp(f)\right)=\supp(h^{-1}\circ f\circ h).$$
That is, the supports of two conjugate self-maps are homeomorphic to each other.
\begin{proof}
Since $h$ is a homeomorphism,
$\supp(f)\cong h^{-1}\left(\supp(f)\right)$. We show that $h^{-1}\left(\supp(f)\right)=\supp(h^{-1}\circ f\circ h)$ :
$$x\in h^{-1}\left(\supp(f)\right)\Leftrightarrow h(x)\in\supp(f)  \Leftrightarrow f(h(x)) \neq h(x) \Leftrightarrow  x\in \supp(h^{-1}\circ f\circ h).$$
\end{proof}
\end{lemma}

Take a representative $\tau_1$ of $\phi(\beta_1)$ and $\widetilde\delta$ of $\phi(\delta)$, where $\delta=\beta_1\cdots\beta_{n-1}$.
Define $\tau_2,\ldots,\tau_{n-1}$ to be $\tau_i:=\widetilde\delta^{-1}\circ\tau_{i-1}\circ\widetilde\delta$ inductively,
then we have $\phi(\beta_i)=[\tau_i]$ for $i=1,\ldots,n-1$.

Each $T_i=\overline{\supp(\tau_i)}$ is homeomorphic to a compact surface
whose boundary is fixed pointwise by $\tau_i$, while $\supp(\tau_i)$ is an open subset of $S_{g,k}$.
Because our goal is to show homology triviality of $\phi$ in the stable range,
we may assume that each $T_i$ is connected.
For if $T_i$ consists of two or more components,
considering the Harer-Ivanov stability theorem,
we may take only one component and push another ones far away to the other side of the big surface.
The self-homeomorphism $\tau_i$ is also assumed to be a simple twist which is defined as follows.

\begin{definition}
\label{def:half_twist}
A self-homeomorphism $\tau_i$ of a connected surface $T_i$ interchanging two identical subsurfaces (or points) $I_i$ and $I_{i+1}$ is called a {\em simple twist} if
\begin{enumerate}[(i)]
    \item $T_i=\overline{\supp(\tau_i)}$ while the boundary of $T_i$ is fixed pointwise by $\tau_i$,
    \item $\tau_i$ is not isotopic to a composition of two nontrivial self-homeomorphisms $f,g$ with $\supp(f)\cap\supp(g)=\emptyset$, that is, $[\tau_i]=[f\circ g]=[g\circ f]$ implies that either $[f]=1$ or $[g]=1$,
    \item if $[\tau_i]=[\gamma_i]^k$ for some self-homeomorphism $\gamma_i$, then $k=1$ or $-1$.
\end{enumerate}
\end{definition}

\begin{figure}[h]
    \centering
    \def \globalscale {1.000000}
    \begin{tikzpicture}[y=0.80pt, x=0.80pt, yscale=-\globalscale, xscale=\globalscale, inner sep=0pt, outer sep=0pt]
      \path[draw=black,line width=0.800pt] (70.0533,79.1502) .. controls (73.1476,82.5177) and (75.6000,86.0136) .. (75.6000,90.0000);

      \path[draw=black,line width=0.800pt] (57.6000,61.2000) .. controls (57.6000,65.1647) and (60.0257,68.6442) .. (63.0962,71.9948);

      \path[draw=black,line width=0.800pt] (75.6000,61.2000) .. controls (75.6000,72.0000) and (57.6000,79.2000) .. (57.6000,90.0000);

      \path[draw=black,line width=0.800pt] (93.6000,61.2000) -- (93.6000,90.0000);

      \path[draw=black,line width=0.800pt] (88.0533,107.9500) .. controls (91.1476,111.3180) and (93.6000,114.8140) .. (93.6000,118.8000);

      \path[draw=black,line width=0.800pt] (75.6000,90.0000) .. controls (75.6000,93.9647) and (78.0257,97.4442) .. (81.0962,100.7950);

      \path[draw=black,line width=0.800pt] (93.6000,90.0000) .. controls (93.6000,100.8000) and (75.6000,108.0000) .. (75.6000,118.8000);

      \path[draw=black,line width=0.800pt] (57.6000,90.0000) -- (57.6000,118.8000);

      \path[draw=black,line width=0.800pt] (70.0533,136.7500) .. controls (73.1476,140.1180) and (75.6000,143.6140) .. (75.6000,147.6000);

      \path[draw=black,line width=0.800pt] (57.6000,118.8000) .. controls (57.6000,122.7650) and (60.0257,126.2440) .. (63.0962,129.5950);

      \path[draw=black,line width=0.800pt] (75.6000,118.8000) .. controls (75.6000,129.6000) and (57.6000,136.8000) .. (57.6000,147.6000);

      \path[draw=black,line width=0.800pt] (93.6000,118.8000) -- (93.6000,147.6000);

      \path[draw=black,line width=0.800pt] (156.4530,79.1502) .. controls (159.5480,82.5177) and (162.0000,86.0136) .. (162.0000,90.0000);

      \path[draw=black,line width=0.800pt] (144.0000,61.2000) .. controls (144.0000,65.1647) and (146.4260,68.6442) .. (149.4960,71.9948);

      \path[draw=black,line width=0.800pt] (162.0000,61.2000) .. controls (162.0000,72.0000) and (144.0000,79.2000) .. (144.0000,90.0000);

      \path[draw=black,line width=0.800pt] (126.0000,61.2000) -- (126.0000,90.0000);

      \path[draw=black,line width=0.800pt] (138.4530,107.9500) .. controls (141.5480,111.3180) and (144.0000,114.8140) .. (144.0000,118.8000);

      \path[draw=black,line width=0.800pt] (126.0000,90.0000) .. controls (126.0000,93.9647) and (128.4260,97.4442) .. (131.4960,100.7950);

      \path[draw=black,line width=0.800pt] (144.0000,90.0000) .. controls (144.0000,100.8000) and (126.0000,108.0000) .. (126.0000,118.8000);

      \path[draw=black,line width=0.800pt] (162.0000,90.0000) -- (162.0000,118.8000);

      \path[draw=black,line width=0.800pt] (156.4530,136.7500) .. controls (159.5480,140.1180) and (162.0000,143.6140) .. (162.0000,147.6000);

      \path[draw=black,line width=0.800pt] (144.0000,118.8000) .. controls (144.0000,122.7650) and (146.4260,126.2440) .. (149.4960,129.5950);

      \path[draw=black,line width=0.800pt] (162.0000,118.8000) .. controls (162.0000,129.6000) and (144.0000,136.8000) .. (144.0000,147.6000);

      \path[draw=black,line width=0.800pt] (126.0000,118.8000) -- (126.0000,147.6000);

      \path[cm={{1.0,0.0,0.0,1.0,(104.4,104.4)}},fill=black] (0.0000,0.0000) node[below right] () {$\cong$};

      \path[draw=black,line width=0.800pt] (237.6010,63.7246) .. controls (241.4270,65.2831) and (244.8000,66.8415) .. (244.8000,68.4000);

      \path[draw=black,line width=0.800pt] (226.8000,57.6000) .. controls (226.8000,59.1583) and (230.1730,60.7166) .. (233.9980,62.2749);

      \path[draw=black,line width=0.800pt] (244.8000,57.6000) .. controls (244.8000,61.2000) and (226.8000,64.8000) .. (226.8000,68.4000);

      \path[draw=black,line width=0.800pt] (262.8000,57.6000) -- (262.8000,68.4000);

      \path[draw=black,line width=0.800pt] (237.6010,74.5246) .. controls (241.4270,76.0831) and (244.8000,77.6415) .. (244.8000,79.2000);

      \path[draw=black,line width=0.800pt] (226.8000,68.4000) .. controls (226.8000,69.9583) and (230.1730,71.5166) .. (233.9980,73.0750);

      \path[draw=black,line width=0.800pt] (244.8000,68.4000) .. controls (244.8000,72.0000) and (226.8000,75.6000) .. (226.8000,79.2000);

      \path[draw=black,line width=0.800pt] (262.8000,68.4000) -- (262.8000,79.2000);

      \path[draw=black,line width=0.800pt] (237.6010,85.3246) .. controls (241.4270,86.8831) and (244.8000,88.4415) .. (244.8000,90.0000);

      \path[draw=black,line width=0.800pt] (226.8000,79.2000) .. controls (226.8000,80.7583) and (230.1730,82.3166) .. (233.9980,83.8750);

      \path[draw=black,line width=0.800pt] (244.8000,79.2000) .. controls (244.8000,82.8000) and (226.8000,86.4000) .. (226.8000,90.0000);

      \path[draw=black,line width=0.800pt] (262.8000,79.2000) -- (262.8000,90.0000);

      \path[draw=black,line width=0.800pt] (255.6010,96.1246) .. controls (259.4270,97.6831) and (262.8000,99.2415) .. (262.8000,100.8000);

      \path[draw=black,line width=0.800pt] (244.8000,90.0000) .. controls (244.8000,91.5583) and (248.1730,93.1166) .. (251.9980,94.6749);

      \path[draw=black,line width=0.800pt] (262.8000,90.0000) .. controls (262.8000,93.6000) and (244.8000,97.2000) .. (244.8000,100.8000);

      \path[draw=black,line width=0.800pt] (226.8000,90.0000) -- (226.8000,100.8000);

      \path[draw=black,line width=0.800pt] (255.6010,106.9250) .. controls (259.4270,108.4830) and (262.8000,110.0420) .. (262.8000,111.6000);

      \path[draw=black,line width=0.800pt] (244.8000,100.8000) .. controls (244.8000,102.3580) and (248.1730,103.9170) .. (251.9980,105.4750);

      \path[draw=black,line width=0.800pt] (262.8000,100.8000) .. controls (262.8000,104.4000) and (244.8000,108.0000) .. (244.8000,111.6000);

      \path[draw=black,line width=0.800pt] (226.8000,100.8000) -- (226.8000,111.6000);

      \path[draw=black,line width=0.800pt] (255.6010,117.7250) .. controls (259.4270,119.2830) and (262.8000,120.8420) .. (262.8000,122.4000);

      \path[draw=black,line width=0.800pt] (244.8000,111.6000) .. controls (244.8000,113.1580) and (248.1730,114.7170) .. (251.9980,116.2750);

      \path[draw=black,line width=0.800pt] (262.8000,111.6000) .. controls (262.8000,115.2000) and (244.8000,118.8000) .. (244.8000,122.4000);

      \path[draw=black,line width=0.800pt] (226.8000,111.6000) -- (226.8000,122.4000);

      \path[draw=black,line width=0.800pt] (237.6010,128.5250) .. controls (241.4270,130.0830) and (244.8000,131.6420) .. (244.8000,133.2000);

      \path[draw=black,line width=0.800pt] (226.8000,122.4000) .. controls (226.8000,123.9580) and (230.1730,125.5170) .. (233.9980,127.0750);

      \path[draw=black,line width=0.800pt] (244.8000,122.4000) .. controls (244.8000,126.0000) and (226.8000,129.6000) .. (226.8000,133.2000);

      \path[draw=black,line width=0.800pt] (262.8000,122.4000) -- (262.8000,133.2000);

      \path[draw=black,line width=0.800pt] (237.6010,139.3250) .. controls (241.4270,140.8830) and (244.8000,142.4420) .. (244.8000,144.0000);

      \path[draw=black,line width=0.800pt] (226.8000,133.2000) .. controls (226.8000,134.7580) and (230.1730,136.3170) .. (233.9980,137.8750);

      \path[draw=black,line width=0.800pt] (244.8000,133.2000) .. controls (244.8000,136.8000) and (226.8000,140.4000) .. (226.8000,144.0000);

      \path[draw=black,line width=0.800pt] (262.8000,133.2000) -- (262.8000,144.0000);

      \path[draw=black,line width=0.800pt] (237.6010,150.1250) .. controls (241.4270,151.6830) and (244.8000,153.2420) .. (244.8000,154.8000);

      \path[draw=black,line width=0.800pt] (226.8000,144.0000) .. controls (226.8000,145.5580) and (230.1730,147.1170) .. (233.9980,148.6750);

      \path[draw=black,line width=0.800pt] (244.8000,144.0000) .. controls (244.8000,147.6000) and (226.8000,151.2000) .. (226.8000,154.8000);

      \path[draw=black,line width=0.800pt] (262.8000,144.0000) -- (262.8000,154.8000);

      \path[draw=black,line width=0.800pt] (324.0010,63.7246) .. controls (327.8270,65.2831) and (331.2000,66.8415) .. (331.2000,68.4000);

      \path[draw=black,line width=0.800pt] (313.2000,57.6000) .. controls (313.2000,59.1583) and (316.5730,60.7166) .. (320.3980,62.2749);

      \path[draw=black,line width=0.800pt] (331.2000,57.6000) .. controls (331.2000,61.2000) and (313.2000,64.8000) .. (313.2000,68.4000);

      \path[draw=black,line width=0.800pt] (295.2000,57.6000) -- (295.2000,68.4000);

      \path[draw=black,line width=0.800pt] (324.0010,74.5246) .. controls (327.8270,76.0831) and (331.2000,77.6415) .. (331.2000,79.2000);

      \path[draw=black,line width=0.800pt] (313.2000,68.4000) .. controls (313.2000,69.9583) and (316.5730,71.5166) .. (320.3980,73.0749);

      \path[draw=black,line width=0.800pt] (331.2000,68.4000) .. controls (331.2000,72.0000) and (313.2000,75.6000) .. (313.2000,79.2000);

      \path[draw=black,line width=0.800pt] (295.2000,68.4000) -- (295.2000,79.2000);

      \path[draw=black,line width=0.800pt] (324.0010,85.3246) .. controls (327.8270,86.8831) and (331.2000,88.4415) .. (331.2000,90.0000);

      \path[draw=black,line width=0.800pt] (313.2000,79.2000) .. controls (313.2000,80.7583) and (316.5730,82.3166) .. (320.3980,83.8749);

      \path[draw=black,line width=0.800pt] (331.2000,79.2000) .. controls (331.2000,82.8000) and (313.2000,86.4000) .. (313.2000,90.0000);

      \path[draw=black,line width=0.800pt] (295.2000,79.2000) -- (295.2000,90.0000);

      \path[draw=black,line width=0.800pt] (306.0010,96.1246) .. controls (309.8270,97.6831) and (313.2000,99.2415) .. (313.2000,100.8000);

      \path[draw=black,line width=0.800pt] (295.2000,90.0000) .. controls (295.2000,91.5583) and (298.5730,93.1166) .. (302.3980,94.6749);

      \path[draw=black,line width=0.800pt] (313.2000,90.0000) .. controls (313.2000,93.6000) and (295.2000,97.2000) .. (295.2000,100.8000);

      \path[draw=black,line width=0.800pt] (331.2000,90.0000) -- (331.2000,100.8000);

      \path[draw=black,line width=0.800pt] (306.0010,106.9250) .. controls (309.8270,108.4830) and (313.2000,110.0420) .. (313.2000,111.6000);

      \path[draw=black,line width=0.800pt] (295.2000,100.8000) .. controls (295.2000,102.3580) and (298.5730,103.9170) .. (302.3980,105.4750);

      \path[draw=black,line width=0.800pt] (313.2000,100.8000) .. controls (313.2000,104.4000) and (295.2000,108.0000) .. (295.2000,111.6000);

      \path[draw=black,line width=0.800pt] (331.2000,100.8000) -- (331.2000,111.6000);

      \path[draw=black,line width=0.800pt] (306.0010,117.7250) .. controls (309.8270,119.2830) and (313.2000,120.8420) .. (313.2000,122.4000);

      \path[draw=black,line width=0.800pt] (295.2000,111.6000) .. controls (295.2000,113.1580) and (298.5730,114.7170) .. (302.3980,116.2750);

      \path[draw=black,line width=0.800pt] (313.2000,111.6000) .. controls (313.2000,115.2000) and (295.2000,118.8000) .. (295.2000,122.4000);

      \path[draw=black,line width=0.800pt] (331.2000,111.6000) -- (331.2000,122.4000);

      \path[draw=black,line width=0.800pt] (324.0010,128.5250) .. controls (327.8270,130.0830) and (331.2000,131.6420) .. (331.2000,133.2000);

      \path[draw=black,line width=0.800pt] (313.2000,122.4000) .. controls (313.2000,123.9580) and (316.5730,125.5170) .. (320.3980,127.0750);

      \path[draw=black,line width=0.800pt] (331.2000,122.4000) .. controls (331.2000,126.0000) and (313.2000,129.6000) .. (313.2000,133.2000);

      \path[draw=black,line width=0.800pt] (295.2000,122.4000) -- (295.2000,133.2000);

      \path[draw=black,line width=0.800pt] (324.0010,139.3250) .. controls (327.8270,140.8830) and (331.2000,142.4420) .. (331.2000,144.0000);

      \path[draw=black,line width=0.800pt] (313.2000,133.2000) .. controls (313.2000,134.7580) and (316.5730,136.3170) .. (320.3980,137.8750);

      \path[draw=black,line width=0.800pt] (331.2000,133.2000) .. controls (331.2000,136.8000) and (313.2000,140.4000) .. (313.2000,144.0000);

      \path[draw=black,line width=0.800pt] (295.2000,133.2000) -- (295.2000,144.0000);

      \path[draw=black,line width=0.800pt] (324.0010,150.1250) .. controls (327.8270,151.6830) and (331.2000,153.2420) .. (331.2000,154.8000);

      \path[draw=black,line width=0.800pt] (313.2000,144.0000) .. controls (313.2000,145.5580) and (316.5730,147.1170) .. (320.3980,148.6750);

      \path[draw=black,line width=0.800pt] (331.2000,144.0000) .. controls (331.2000,147.6000) and (313.2000,151.2000) .. (313.2000,154.8000);

      \path[draw=black,line width=0.800pt] (295.2000,144.0000) -- (295.2000,154.8000);

      \path[cm={{1.0,0.0,0.0,1.0,(273.6,104.4)}},fill=black] (0.0000,0.0000) node[below right] () {$\not\cong$};

    \end{tikzpicture}
    \caption{Simple twists $\tau_1,\tau_2$ satisfy the braid relation, but $\tau_1^3,\tau_2^3$ do not.}
    \label{fig:simple}
\end{figure}
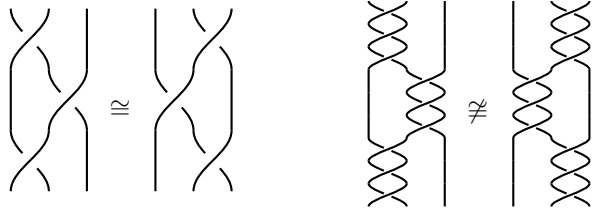

The condition (ii) implies the minimality of the choice of the atomic surface $T_i$.
The condition (iii) is necessary because for a pair of standard generators $\beta_1,\beta_2$ of $B_n$ satisfy the braid relation $\beta_1\beta_2\beta_1=\beta_2\beta_1\beta_2$, but their odd products $\beta_1^{2i+1},\beta_2^{2i+1}$ do not satisfy the braid relation although they interchange two points by a twist (see Figure~\ref{fig:simple}).

\begin{definition}
\label{def:regular}
An embedding $\phi:B_n\hookrightarrow\Gamma_{g,k}, \beta_i\mapsto[\tau_i]$ is a {\em regular embedding} if
\begin{enumerate}[(i)]
    \item $T_i=\overline{\supp(\tau_i)}$ are all identical and connected,
    \item each $\tau_i$ is a simple twist interchanging two identical subsurfaces (or points) $I_i$ and $I_{i+1}$ of $T_i$,
    \item $I_i\subset T_{i-1}\cap T_i$ for $i=2,\ldots,n$ and $S_{g,k}\cong\cup_{i=1}^{n-1}T_i$.
\end{enumerate}
\end{definition}

Here we give three examples of regular (non-geometric) embeddings.

\begin{example}
\label{ex:dfold}
As seen in the section~\ref{sec:d-fold}, there is an infinite family of non-geometric embeddings
induced by the lifts of half Dehn twists on disk with some marked points
via $d$-fold covering map.
Let $p:S_{g,b}^{(n)}\rightarrow S_{0,1}^{(n)}$ be a cyclic $d$-fold branched covering
with $n$ branch points.
Then an embedding $\phi:B_n\hookrightarrow\Gamma_{g,b},\beta_i\mapsto[\tau_i]$ is induced from this $d$-fold covering;
where $\tau_i$ is obtained by a lift of half Dehn twist on disk with two marked points
via covering map $p$:
for $p^{-1}(S_{0,1}^{(2)})=S_{h,k}^{(2)}$ we have
$$
\xymatrix{
S_{h,k}^{(2)} \ar[d]_p \ar[r]^{\tau_i} & S_{h,k}^{(2)}\ar[d]^p\\
S_{0,1}^{(2)} \ar[r]^{\beta_i} & S_{0,1}^{(2)}
}
$$

For $d=2$, this embedding is the usual geometric embedding.
In the case $d\geq 3$, these embeddings are all non-geometric.
These embeddings are all regular embeddings with the atomic surfaces $T_i\cong p^{-1}(S_{0,1}^{(2)})=S_{h,k}^{(2)}$,
and each $I_i$ is a point.
\end{example}

\begin{example}[Szepietowski's construction \cite{Szepietowski10}]
Szepietowski constructed a non-geometric embedding $B_n\hookrightarrow\Gamma_{n-1,2}$
using non-orientable surface $N_{n,1}$
and the 2-fold covering map $p:S_{n-1,2}\rightarrow N_{n,1}$
induced by mapping two antipodal points of $S_{1,2}$ to one point of Möbius band.
Starting with $S_{0,(n)+1}$, the disk with $n$ holes,
we paste to each hole a Möbius band $N_{1,1}$ and get a non-orientable surface $N_{n,1}$.
Then we get a natural embedding $$\varphi:B_n\hookrightarrow\mathcal N_{n,1}$$
where $\mathcal N_{n,1}$ denotes the mapping class group of $N_{n,1}$.
The lift of homeomorphisms of $N_{n,1}$ to its double covering space $S_{n-1,2}$
induces an injection $L:\mathcal N_{n,1}\hookrightarrow\Gamma_{n-1,2}$.
Then the composite $L\circ\varphi:B_n\hookrightarrow\Gamma_{n-1,2}$ is a non-geometric regular embedding.
The injectivity of the map $L$ can be shown by the result of Birman and Chillingworth (\cite{BC72}).
\end{example}

\begin{example}[The mirror construction \cite{ST07},\cite{BT12}]
The obvious injection $B_n=\Gamma_{0,1}^{(n)}\hookrightarrow\Gamma_{0,(n)+1}$ can be extended to an injection $\phi:B_n\hookrightarrow\Gamma_{n-1,2}$ by the mirror construction.
From two sheets of disks with $n$ holes,
by gluing them along the boundaries of the holes,
we get a surface $S_{n-1,2}$ (see Figure~\ref{fig:mirror}).

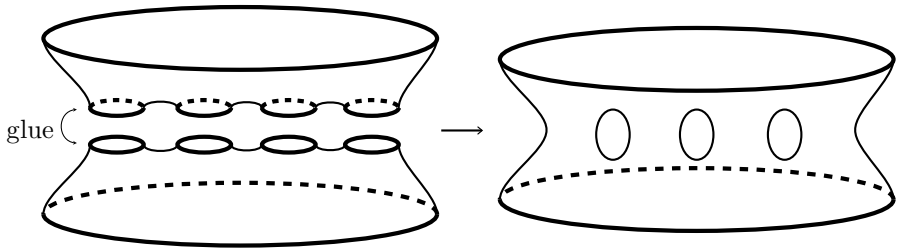
\begin{figure}[h]
    \centering
    \def \globalscale {0.55}
    \begin{tikzpicture}[y=0.80pt, x=0.80pt, yscale=-\globalscale, xscale=\globalscale, inner sep=0pt, outer sep=0pt]
      \path[draw=black,line width=1.600pt] (21.6000,219.6000) .. controls (21.6000,204.6880) and (97.3535,192.6000) .. (190.8000,192.6000) .. controls (284.2470,192.6000) and (360.0000,204.6880) .. (360.0000,219.6000) .. controls (360.0000,234.5120) and (284.2470,246.6000) .. (190.8000,246.6000) .. controls (97.3535,246.6000) and (21.6000,234.5120) .. (21.6000,219.6000) -- cycle;

      \path[draw=black,line width=0.800pt] (21.6000,219.6000) .. controls (21.6000,237.6030) and (62.1120,261.4210) .. (61.6930,279.4200);

      \path[draw=black,line width=0.800pt] (107.8000,279.4200) .. controls (107.7180,272.1700) and (136.8000,272.1700) .. (136.8000,279.4200);

      \path[draw=black,line width=0.800pt] (182.9070,279.4200) .. controls (182.7210,272.7760) and (209.7100,272.7770) .. (209.4930,279.4200);

      \path[draw=black,line width=0.800pt] (255.6000,279.4200) .. controls (255.6000,273.0830) and (280.9990,273.0830) .. (280.9470,279.4200);

      \path[draw=black,line width=0.800pt] (360.0000,219.6000) .. controls (360.0000,236.6730) and (326.6130,262.3520) .. (327.0530,279.4200);

      \path[draw=black,line width=1.600pt] (61.6930,310.9800) .. controls (61.6930,314.7250) and (72.0144,317.7610) .. (84.7465,317.7610) .. controls (97.4786,317.7610) and (107.8000,314.7250) .. (107.8000,310.9800) .. controls (107.8000,307.2360) and (97.4786,304.2000) .. (84.7465,304.2000) .. controls (72.0144,304.2000) and (61.6930,307.2360) .. (61.6930,310.9800) -- cycle;

      \path[draw=black,line width=0.800pt] (21.6000,370.8000) .. controls (21.6000,352.7970) and (62.1120,328.9790) .. (61.6930,310.9800);

      \path[draw=black,line width=1.600pt] (136.8000,310.9800) .. controls (136.8000,314.7250) and (147.1210,317.7610) .. (159.8540,317.7610) .. controls (172.5860,317.7610) and (182.9070,314.7250) .. (182.9070,310.9800) .. controls (182.9070,307.2360) and (172.5860,304.2000) .. (159.8540,304.2000) .. controls (147.1210,304.2000) and (136.8000,307.2360) .. (136.8000,310.9800) -- cycle;

      \path[draw=black,line width=1.600pt] (280.9470,310.9800) .. controls (280.9470,314.7250) and (291.2680,317.7610) .. (304.0000,317.7610) .. controls (316.7320,317.7610) and (327.0530,314.7250) .. (327.0530,310.9800) .. controls (327.0530,307.2360) and (316.7320,304.2000) .. (304.0000,304.2000) .. controls (291.2680,304.2000) and (280.9470,307.2360) .. (280.9470,310.9800) -- cycle;

      \path[draw=black,line width=1.600pt] (209.4930,310.9800) .. controls (209.4930,314.7250) and (219.8140,317.7610) .. (232.5470,317.7610) .. controls (245.2790,317.7610) and (255.6000,314.7250) .. (255.6000,310.9800) .. controls (255.6000,307.2360) and (245.2790,304.2000) .. (232.5470,304.2000) .. controls (219.8140,304.2000) and (209.4930,307.2360) .. (209.4930,310.9800) -- cycle;

      \path[draw=black,line width=0.800pt] (107.8000,310.9800) .. controls (107.7180,318.2300) and (136.8000,318.2300) .. (136.8000,310.9800);

      \path[draw=black,line width=0.800pt] (182.9070,310.9800) .. controls (182.7210,317.6240) and (209.7100,317.6230) .. (209.4930,310.9800);

      \path[draw=black,line width=0.800pt] (255.6000,310.9800) .. controls (255.6000,317.3170) and (280.9990,317.3170) .. (280.9470,310.9800);

      \path[draw=black,line width=0.800pt] (360.0000,370.8000) .. controls (360.0000,353.7270) and (326.6130,328.0480) .. (327.0530,310.9800);

      \path[draw=black,dash pattern=on 3.20pt,line width=1.600pt] (360.0000,370.8000) .. controls (360.0000,355.8880) and (284.2470,343.8000) .. (190.8000,343.8000) .. controls (97.3535,343.8000) and (21.6000,355.8880) .. (21.6000,370.8000);

      \path[draw=black,line width=1.600pt] (21.6000,370.8000) .. controls (21.6000,385.7120) and (97.3535,397.8000) .. (190.8000,397.8000) .. controls (284.2470,397.8000) and (360.0000,385.7120) .. (360.0000,370.8000);

      \path[draw=black,line width=1.600pt] (107.8000,279.4200) .. controls (107.8000,283.1640) and (97.4786,286.2000) .. (84.7465,286.2000) .. controls (72.0144,286.2000) and (61.6930,283.1640) .. (61.6930,279.4200);

      \path[draw=black,dash pattern=on 2.40pt,line width=1.600pt] (61.6930,279.4200) .. controls (61.6930,275.6750) and (72.0144,272.6390) .. (84.7465,272.6390) .. controls (97.4786,272.6390) and (107.8000,275.6750) .. (107.8000,279.4200);

      \path[draw=black,line width=1.600pt] (182.9070,279.4200) .. controls (182.9070,283.1640) and (172.5860,286.2000) .. (159.8540,286.2000) .. controls (147.1210,286.2000) and (136.8000,283.1640) .. (136.8000,279.4200);

      \path[draw=black,dash pattern=on 2.40pt,line width=1.600pt] (136.8000,279.4200) .. controls (136.8000,275.6750) and (147.1210,272.6390) .. (159.8540,272.6390) .. controls (172.5860,272.6390) and (182.9070,275.6750) .. (182.9070,279.4200);

      \path[draw=black,line width=1.600pt] (327.0530,279.4200) .. controls (327.0530,283.1640) and (316.7320,286.2000) .. (304.0000,286.2000) .. controls (291.2680,286.2000) and (280.9470,283.1640) .. (280.9470,279.4200);

      \path[draw=black,dash pattern=on 2.40pt,line width=1.600pt] (280.9470,279.4200) .. controls (280.9470,275.6750) and (291.2680,272.6390) .. (304.0000,272.6390) .. controls (316.7320,272.6390) and (327.0530,275.6750) .. (327.0530,279.4200);

      \path[draw=black,line width=1.600pt] (255.6000,279.4200) .. controls (255.6000,283.1640) and (245.2790,286.2000) .. (232.5470,286.2000) .. controls (219.8140,286.2000) and (209.4930,283.1640) .. (209.4930,279.4200);

      \path[draw=black,dash pattern=on 2.40pt,line width=1.600pt] (209.4930,279.4200) .. controls (209.4930,275.6750) and (219.8140,272.6390) .. (232.5470,272.6390) .. controls (245.2790,272.6390) and (255.6000,275.6750) .. (255.6000,279.4200);

      \path[draw=black,line width=1.600pt] (414.0000,237.6000) .. controls (414.0000,222.6880) and (489.7530,210.6000) .. (583.2000,210.6000) .. controls (676.6470,210.6000) and (752.4000,222.6880) .. (752.4000,237.6000) .. controls (752.4000,252.5120) and (676.6470,264.6000) .. (583.2000,264.6000) .. controls (489.7530,264.6000) and (414.0000,252.5120) .. (414.0000,237.6000) -- cycle;

      \path[draw=black,line width=0.800pt] (414.0000,237.6000) .. controls (414.0000,255.6030) and (454.5120,280.3820) .. (454.0930,298.3800);

      \path[draw=black,line width=0.800pt] (752.4000,237.6000) .. controls (752.4000,254.6730) and (719.0130,281.3130) .. (719.4530,298.3800);

      \path[draw=black,line width=0.800pt] (414.0000,358.2000) .. controls (414.0000,340.1970) and (454.5120,316.3790) .. (454.0930,298.3800);

      \path[draw=black,line width=0.800pt] (752.4000,358.2000) .. controls (752.4000,341.1270) and (719.0130,315.4480) .. (719.4530,298.3800);

      \path[draw=black,dash pattern=on 3.20pt,line width=1.600pt] (752.4000,358.2000) .. controls (752.4000,343.2880) and (676.6470,331.2000) .. (583.2000,331.2000) .. controls (489.7530,331.2000) and (414.0000,343.2880) .. (414.0000,358.2000);

      \path[draw=black,line width=1.600pt] (414.0000,358.2000) .. controls (414.0000,373.1120) and (489.7530,385.2000) .. (583.2000,385.2000) .. controls (676.6470,385.2000) and (752.4000,373.1120) .. (752.4000,358.2000);

      \path[draw=black,line width=0.800pt] (496.8000,302.4000) .. controls (496.8000,290.4710) and (503.2470,280.8000) .. (511.2000,280.8000) .. controls (519.1530,280.8000) and (525.6000,290.4710) .. (525.6000,302.4000) .. controls (525.6000,314.3290) and (519.1530,324.0000) .. (511.2000,324.0000) .. controls (503.2470,324.0000) and (496.8000,314.3290) .. (496.8000,302.4000) -- cycle;

      \path[draw=black,line width=0.800pt] (568.8000,302.4000) .. controls (568.8000,290.4710) and (575.2470,280.8000) .. (583.2000,280.8000) .. controls (591.1530,280.8000) and (597.6000,290.4710) .. (597.6000,302.4000) .. controls (597.6000,314.3290) and (591.1530,324.0000) .. (583.2000,324.0000) .. controls (575.2470,324.0000) and (568.8000,314.3290) .. (568.8000,302.4000) -- cycle;

      \path[draw=black,line width=0.800pt] (644.4000,302.4000) .. controls (644.4000,290.4710) and (650.8470,280.8000) .. (658.8000,280.8000) .. controls (666.7530,280.8000) and (673.2000,290.4710) .. (673.2000,302.4000) .. controls (673.2000,314.3290) and (666.7530,324.0000) .. (658.8000,324.0000) .. controls (650.8470,324.0000) and (644.4000,314.3290) .. (644.4000,302.4000) -- cycle;

          \path[draw=black,line width=0.400pt] (48.6496,280.8950) .. controls (32.9820,282.6340) and (32.9826,307.7680) .. (48.6516,309.5050);

          \path[fill=black] (48.2147,309.4810) -- (46.8094,311.1580) -- (50.3990,309.6000) -- (46.9989,307.6630) -- (48.2147,309.4810) -- cycle;

          \path[fill=black] (48.2128,280.9190) -- (46.9970,282.7370) -- (50.3970,280.8000) -- (46.8074,279.2420) -- (48.2128,280.9190) -- cycle;

      \path[cm={{1.0,0.0,0.0,1.0,(-10.0,289.7)}},fill=black] (0.0000,0.0000) node[below right] () {glue};

          \path[draw=black,line width=0.800pt] (398.3100,298.8000) .. controls (395.5090,298.8000) and (363.6000,298.8000) .. (363.6000,298.8000);

          \path[fill=black] (395.4600,302.1870) -- (399.1260,299.1870) -- (399.5990,298.8000) -- (399.1260,298.4130) -- (395.4600,295.4130) .. controls (395.2460,295.2380) and (394.9310,295.2700) .. (394.7560,295.4830) .. controls (394.5810,295.6970) and (394.6130,296.0120) .. (394.8260,296.1870) -- (398.4930,299.1870) -- (398.4930,298.4130) -- (394.8260,301.4130) .. controls (394.6130,301.5880) and (394.5810,301.9030) .. (394.7560,302.1170) .. controls (394.9310,302.3300) and (395.2460,302.3620) .. (395.4600,302.1870) -- cycle;

    \end{tikzpicture}
    \caption{Mirror construction}
    \label{fig:mirror}
\end{figure}

This construction gives rise to an embedding $\phi:B_n\hookrightarrow\Gamma_{n-1,2},\beta_i\mapsto\tau_i$,
where $\tau_i$ is called the pillar switching (\cite{JS13}).
$\phi$ is obviously a regular embedding with
the atomic surface $T_i\cong S_{1,2}$ and
the interchangeable subsurface $I_i$ is the cylinder $S_{0,2}$.
This $\phi$ is essentially equal
to that obtained by Szepietowski's construction.

Note that we may have more general mirror constructions
by taking the interchangeable subsurface $I_i$ as any surface of type $S_{g,2}$ (Figure~\ref{fig:general_mirror}).

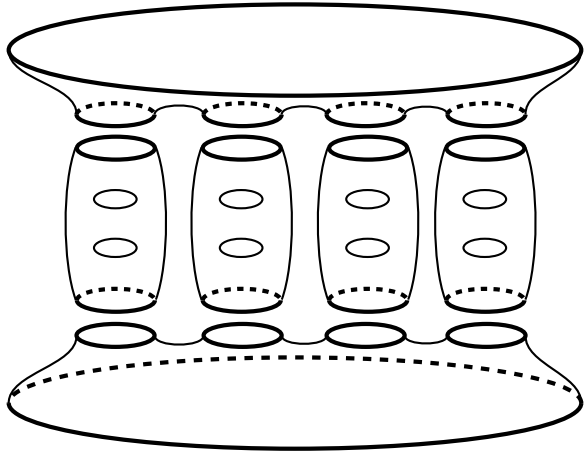
\begin{figure}[h]
    \centering
    \def \globalscale {0.800000}
    \begin{tikzpicture}[y=0.80pt, x=0.80pt, yscale=-\globalscale, xscale=\globalscale, inner sep=0pt, outer sep=0pt]
      \path[draw=black,line width=1.600pt] (21.6000,241.2000) .. controls (21.6000,226.2880) and (97.3534,214.2000) .. (190.8000,214.2000) .. controls (284.2470,214.2000) and (360.0000,226.2880) .. (360.0000,241.2000) .. controls (360.0000,256.1120) and (284.2470,268.2000) .. (190.8000,268.2000) .. controls (97.3534,268.2000) and (21.6000,256.1120) .. (21.6000,241.2000) -- cycle;

      \path[draw=black,line width=0.800pt] (21.6000,241.2000) .. controls (21.6000,259.2030) and (62.1120,261.4210) .. (61.6930,279.4200);

      \path[draw=black,line width=0.800pt] (107.8000,279.4200) .. controls (107.7180,272.1700) and (136.8000,272.1700) .. (136.8000,279.4200);

      \path[draw=black,line width=0.800pt] (182.9070,279.4200) .. controls (182.7210,272.7760) and (209.7100,272.7770) .. (209.4930,279.4200);

      \path[draw=black,line width=0.800pt] (255.6000,279.4200) .. controls (255.6000,273.0830) and (280.9990,273.0830) .. (280.9470,279.4200);

      \path[draw=black,line width=0.800pt] (360.0000,241.2000) .. controls (360.0000,258.2730) and (326.6130,262.3520) .. (327.0530,279.4200);

      \path[draw=black,line width=1.600pt] (61.6930,409.9800) .. controls (61.6930,413.7250) and (72.0144,416.7610) .. (84.7465,416.7610) .. controls (97.4786,416.7610) and (107.8000,413.7250) .. (107.8000,409.9800) .. controls (107.8000,406.2360) and (97.4786,403.2000) .. (84.7465,403.2000) .. controls (72.0144,403.2000) and (61.6930,406.2360) .. (61.6930,409.9800) -- cycle;

      \path[draw=black,line width=0.800pt] (21.6000,450.0000) .. controls (21.6000,431.9970) and (62.1120,427.9790) .. (61.6930,409.9800);

      \path[draw=black,line width=1.600pt] (136.8000,409.9800) .. controls (136.8000,413.7250) and (147.1210,416.7610) .. (159.8540,416.7610) .. controls (172.5860,416.7610) and (182.9070,413.7250) .. (182.9070,409.9800) .. controls (182.9070,406.2360) and (172.5860,403.2000) .. (159.8540,403.2000) .. controls (147.1210,403.2000) and (136.8000,406.2360) .. (136.8000,409.9800) -- cycle;

      \path[draw=black,line width=1.600pt] (280.9470,409.9800) .. controls (280.9470,413.7250) and (291.2680,416.7610) .. (304.0000,416.7610) .. controls (316.7320,416.7610) and (327.0540,413.7250) .. (327.0540,409.9800) .. controls (327.0540,406.2360) and (316.7320,403.2000) .. (304.0000,403.2000) .. controls (291.2680,403.2000) and (280.9470,406.2360) .. (280.9470,409.9800) -- cycle;

      \path[draw=black,line width=1.600pt] (209.4930,409.9800) .. controls (209.4930,413.7250) and (219.8140,416.7610) .. (232.5470,416.7610) .. controls (245.2790,416.7610) and (255.6000,413.7250) .. (255.6000,409.9800) .. controls (255.6000,406.2360) and (245.2790,403.2000) .. (232.5470,403.2000) .. controls (219.8140,403.2000) and (209.4930,406.2360) .. (209.4930,409.9800) -- cycle;

      \path[draw=black,line width=0.800pt] (107.8000,409.9800) .. controls (107.7180,417.2300) and (136.8000,417.2300) .. (136.8000,409.9800);

      \path[draw=black,line width=0.800pt] (182.9070,409.9800) .. controls (182.7210,416.6240) and (209.7100,416.6230) .. (209.4930,409.9800);

      \path[draw=black,line width=0.800pt] (255.6000,409.9800) .. controls (255.6000,416.3170) and (280.9990,416.3170) .. (280.9470,409.9800);

      \path[draw=black,line width=0.800pt] (360.0000,450.0000) .. controls (360.0000,432.9270) and (326.6130,427.0480) .. (327.0540,409.9800);

      \path[draw=black,dash pattern=on 3.20pt,line width=1.600pt] (360.0000,450.0000) .. controls (360.0000,435.0880) and (284.2470,423.0000) .. (190.8000,423.0000) .. controls (97.3534,423.0000) and (21.6000,435.0880) .. (21.6000,450.0000);

      \path[draw=black,line width=1.600pt] (21.6000,450.0000) .. controls (21.6000,464.9120) and (97.3534,477.0000) .. (190.8000,477.0000) .. controls (284.2470,477.0000) and (360.0000,464.9120) .. (360.0000,450.0000);

      \path[draw=black,line width=1.600pt] (107.8000,279.4200) .. controls (107.8000,283.1640) and (97.4786,286.2000) .. (84.7465,286.2000) .. controls (72.0144,286.2000) and (61.6930,283.1640) .. (61.6930,279.4200);

      \path[draw=black,dash pattern=on 2.40pt,line width=1.600pt] (61.6930,279.4200) .. controls (61.6930,275.6750) and (72.0144,272.6390) .. (84.7465,272.6390) .. controls (97.4786,272.6390) and (107.8000,275.6750) .. (107.8000,279.4200);

      \path[draw=black,line width=1.600pt] (182.9070,279.4200) .. controls (182.9070,283.1640) and (172.5860,286.2000) .. (159.8540,286.2000) .. controls (147.1210,286.2000) and (136.8000,283.1640) .. (136.8000,279.4200);

      \path[draw=black,dash pattern=on 2.40pt,line width=1.600pt] (136.8000,279.4200) .. controls (136.8000,275.6750) and (147.1210,272.6390) .. (159.8540,272.6390) .. controls (172.5860,272.6390) and (182.9070,275.6750) .. (182.9070,279.4200);

      \path[draw=black,line width=1.600pt] (327.0530,279.4200) .. controls (327.0530,283.1640) and (316.7320,286.2000) .. (304.0000,286.2000) .. controls (291.2680,286.2000) and (280.9470,283.1640) .. (280.9470,279.4200);

      \path[draw=black,dash pattern=on 2.40pt,line width=1.600pt] (280.9470,279.4200) .. controls (280.9470,275.6750) and (291.2680,272.6390) .. (304.0000,272.6390) .. controls (316.7320,272.6390) and (327.0530,275.6750) .. (327.0530,279.4200);

      \path[draw=black,line width=1.600pt] (255.6000,279.4200) .. controls (255.6000,283.1640) and (245.2790,286.2000) .. (232.5470,286.2000) .. controls (219.8140,286.2000) and (209.4930,283.1640) .. (209.4930,279.4200);

      \path[draw=black,dash pattern=on 2.40pt,line width=1.600pt] (209.4930,279.4200) .. controls (209.4930,275.6750) and (219.8140,272.6390) .. (232.5470,272.6390) .. controls (245.2790,272.6390) and (255.6000,275.6750) .. (255.6000,279.4200);

      \path[draw=black,line width=1.600pt] (61.8930,299.2200) .. controls (61.8930,302.9640) and (72.2144,306.0000) .. (84.9465,306.0000) .. controls (97.6786,306.0000) and (108.0000,302.9640) .. (108.0000,299.2200) .. controls (108.0000,295.4750) and (97.6786,292.4390) .. (84.9465,292.4390) .. controls (72.2144,292.4390) and (61.8930,295.4750) .. (61.8930,299.2200) -- cycle;

      \path[draw=black,line width=1.600pt] (108.0000,389.2200) .. controls (108.0000,392.9640) and (97.6786,396.0000) .. (84.9465,396.0000) .. controls (72.2144,396.0000) and (61.8930,392.9640) .. (61.8930,389.2200);

      \path[draw=black,dash pattern=on 2.40pt,line width=1.600pt] (61.8930,389.2200) .. controls (61.8930,385.4750) and (72.2144,382.4390) .. (84.9465,382.4390) .. controls (97.6786,382.4390) and (108.0000,385.4750) .. (108.0000,389.2200);

      \path[draw=black,line width=0.800pt] (61.2000,298.8000) .. controls (53.0384,319.7680) and (53.3992,367.6960) .. (61.2000,388.8000);

      \path[draw=black,line width=0.800pt] (108.6070,298.8000) .. controls (116.7680,319.7680) and (116.4070,367.6960) .. (108.6070,388.8000);

      \path[draw=black,line width=0.800pt] (72.0000,329.4000) .. controls (72.0000,326.4180) and (77.6412,324.0000) .. (84.6000,324.0000) .. controls (91.5588,324.0000) and (97.2000,326.4180) .. (97.2000,329.4000) .. controls (97.2000,332.3820) and (91.5588,334.8000) .. (84.6000,334.8000) .. controls (77.6412,334.8000) and (72.0000,332.3820) .. (72.0000,329.4000) -- cycle;

      \path[draw=black,line width=0.800pt] (72.0000,358.2000) .. controls (72.0000,355.2180) and (77.6412,352.8000) .. (84.6000,352.8000) .. controls (91.5588,352.8000) and (97.2000,355.2180) .. (97.2000,358.2000) .. controls (97.2000,361.1820) and (91.5588,363.6000) .. (84.6000,363.6000) .. controls (77.6412,363.6000) and (72.0000,361.1820) .. (72.0000,358.2000) -- cycle;

      \path[draw=black,line width=1.600pt] (136.2800,299.2200) .. controls (136.2800,302.9640) and (146.6010,306.0000) .. (159.3330,306.0000) .. controls (172.0650,306.0000) and (182.3870,302.9640) .. (182.3870,299.2200) .. controls (182.3870,295.4750) and (172.0650,292.4390) .. (159.3330,292.4390) .. controls (146.6010,292.4390) and (136.2800,295.4750) .. (136.2800,299.2200) -- cycle;

      \path[draw=black,line width=1.600pt] (182.3870,389.2200) .. controls (182.3870,392.9640) and (172.0650,396.0000) .. (159.3330,396.0000) .. controls (146.6010,396.0000) and (136.2800,392.9640) .. (136.2800,389.2200);

      \path[draw=black,dash pattern=on 2.40pt,line width=1.600pt] (136.2800,389.2200) .. controls (136.2800,385.4750) and (146.6010,382.4390) .. (159.3330,382.4390) .. controls (172.0650,382.4390) and (182.3870,385.4750) .. (182.3870,389.2200);

      \path[draw=black,line width=0.800pt] (135.5870,298.8000) .. controls (127.4250,319.7680) and (127.7860,367.6960) .. (135.5870,388.8000);

      \path[draw=black,line width=0.800pt] (182.9930,298.8000) .. controls (191.1550,319.7680) and (190.7940,367.6960) .. (182.9930,388.8000);

      \path[draw=black,line width=0.800pt] (146.3870,329.4000) .. controls (146.3870,326.4180) and (152.0280,324.0000) .. (158.9870,324.0000) .. controls (165.9450,324.0000) and (171.5870,326.4180) .. (171.5870,329.4000) .. controls (171.5870,332.3820) and (165.9450,334.8000) .. (158.9870,334.8000) .. controls (152.0280,334.8000) and (146.3870,332.3820) .. (146.3870,329.4000) -- cycle;

      \path[draw=black,line width=0.800pt] (146.3870,358.2000) .. controls (146.3870,355.2180) and (152.0280,352.8000) .. (158.9870,352.8000) .. controls (165.9450,352.8000) and (171.5870,355.2180) .. (171.5870,358.2000) .. controls (171.5870,361.1820) and (165.9450,363.6000) .. (158.9870,363.6000) .. controls (152.0280,363.6000) and (146.3870,361.1820) .. (146.3870,358.2000) -- cycle;

      \path[draw=black,line width=1.600pt] (210.9900,299.2200) .. controls (210.9900,302.9640) and (221.3110,306.0000) .. (234.0430,306.0000) .. controls (246.7750,306.0000) and (257.0970,302.9640) .. (257.0970,299.2200) .. controls (257.0970,295.4750) and (246.7750,292.4390) .. (234.0430,292.4390) .. controls (221.3110,292.4390) and (210.9900,295.4750) .. (210.9900,299.2200) -- cycle;

      \path[draw=black,line width=1.600pt] (257.0970,389.2200) .. controls (257.0970,392.9640) and (246.7750,396.0000) .. (234.0430,396.0000) .. controls (221.3110,396.0000) and (210.9900,392.9640) .. (210.9900,389.2200);

      \path[draw=black,dash pattern=on 2.40pt,line width=1.600pt] (210.9900,389.2200) .. controls (210.9900,385.4750) and (221.3110,382.4390) .. (234.0430,382.4390) .. controls (246.7750,382.4390) and (257.0970,385.4750) .. (257.0970,389.2200);

      \path[draw=black,line width=0.800pt] (210.2970,298.8000) .. controls (202.1350,319.7680) and (202.4960,367.6960) .. (210.2970,388.8000);

      \path[draw=black,line width=0.800pt] (257.7030,298.8000) .. controls (265.8650,319.7680) and (265.5040,367.6960) .. (257.7030,388.8000);

      \path[draw=black,line width=0.800pt] (221.0970,329.4000) .. controls (221.0970,326.4180) and (226.7380,324.0000) .. (233.6970,324.0000) .. controls (240.6550,324.0000) and (246.2970,326.4180) .. (246.2970,329.4000) .. controls (246.2970,332.3820) and (240.6550,334.8000) .. (233.6970,334.8000) .. controls (226.7380,334.8000) and (221.0970,332.3820) .. (221.0970,329.4000) -- cycle;

      \path[draw=black,line width=0.800pt] (221.0970,358.2000) .. controls (221.0970,355.2180) and (226.7380,352.8000) .. (233.6970,352.8000) .. controls (240.6550,352.8000) and (246.2970,355.2180) .. (246.2970,358.2000) .. controls (246.2970,361.1820) and (240.6550,363.6000) .. (233.6970,363.6000) .. controls (226.7380,363.6000) and (221.0970,361.1820) .. (221.0970,358.2000) -- cycle;

      \path[draw=black,line width=1.600pt] (280.2800,299.2200) .. controls (280.2800,302.9640) and (290.6010,306.0000) .. (303.3330,306.0000) .. controls (316.0650,306.0000) and (326.3870,302.9640) .. (326.3870,299.2200) .. controls (326.3870,295.4750) and (316.0650,292.4390) .. (303.3330,292.4390) .. controls (290.6010,292.4390) and (280.2800,295.4750) .. (280.2800,299.2200) -- cycle;

      \path[draw=black,line width=1.600pt] (326.3870,389.2200) .. controls (326.3870,392.9640) and (316.0650,396.0000) .. (303.3330,396.0000) .. controls (290.6010,396.0000) and (280.2800,392.9640) .. (280.2800,389.2200);

      \path[draw=black,dash pattern=on 2.40pt,line width=1.600pt] (280.2800,389.2200) .. controls (280.2800,385.4750) and (290.6010,382.4390) .. (303.3330,382.4390) .. controls (316.0650,382.4390) and (326.3870,385.4750) .. (326.3870,389.2200);

      \path[draw=black,line width=0.800pt] (279.5870,298.8000) .. controls (271.4250,319.7680) and (271.7860,367.6960) .. (279.5870,388.8000);

      \path[draw=black,line width=0.800pt] (326.9930,298.8000) .. controls (335.1550,319.7680) and (334.7940,367.6960) .. (326.9930,388.8000);

      \path[draw=black,line width=0.800pt] (290.3870,329.4000) .. controls (290.3870,326.4180) and (296.0280,324.0000) .. (302.9870,324.0000) .. controls (309.9450,324.0000) and (315.5870,326.4180) .. (315.5870,329.4000) .. controls (315.5870,332.3820) and (309.9450,334.8000) .. (302.9870,334.8000) .. controls (296.0280,334.8000) and (290.3870,332.3820) .. (290.3870,329.4000) -- cycle;

      \path[draw=black,line width=0.800pt] (290.3870,358.2000) .. controls (290.3870,355.2180) and (296.0280,352.8000) .. (302.9870,352.8000) .. controls (309.9450,352.8000) and (315.5870,355.2180) .. (315.5870,358.2000) .. controls (315.5870,361.1820) and (309.9450,363.6000) .. (302.9870,363.6000) .. controls (296.0280,363.6000) and (290.3870,361.1820) .. (290.3870,358.2000) -- cycle;

    \end{tikzpicture}
    \caption{General mirror construction (Type I).}
    \label{fig:general_mirror}
\end{figure}

More generally, we may have more than two sheets of disks with $n$ holes.
If we have $m$ $(m \ge 3)$ sheets of disks, then we get a mirror construction by gluing $n$ copies of $S_{g,m}$, which are interchangeable surfaces $I_i$, to these disks.
That is, each hole of $S_{g,m}$ is glued to a hole of a disk.
Call this type of mirror construction, Type II.
The embedding induced by the mirror construction of this type is obviously regular. 

There are even more general types of regular embeddings than these mirror constructions.
\begin{itemize}
    \item[III.] On each sheet of disks, the number of holes may not be necessarily all equal.
    \item[IV.] Instead of sheets of disks, we may have any surfaces and any simple twists on them.
\end{itemize}
All of these general regular embeddings will be considered in the next section.
\end{example}

\section{The general Harer conjecture}
\label{sec:Harer}

In consideration of the Harer-Ivanov stability theorem,
it suffices to prove the Harer conjecture for regular embeddings.
We are, as in the proof of Theorem~\ref{thm:harer}, interested in only the number $k$ of boundary components of the surface with respect to a regular embedding $\phi:B_n\hookrightarrow\Gamma_{g,k}, \beta_i\mapsto[\tau_i]$.

\begin{theorem}
\label{thm:GeneralHarer}
For every embedding $\phi : B_{n} \hookrightarrow \Gamma_{g,k}$,
$$\phi_{*}:H_{i}(B_{\infty}; R) \rightarrow H_{i}(\Gamma_{\infty}; R)$$
is trivial for all $i\geq 1$ and any coefficient $R$.
\end{theorem}
\begin{proof}
Let $\phi:B_n\rightarrow\Gamma_{g,k}$ be a regular embedding.
Let us take a look at just $T_1=\overline{\supp(\tau_1)}$ with two interchangeable subsurfaces (or points) $I_1$ and $I_2$. Recall that all $T_i$ are identical for $i=1,\ldots,n-1$.

Let us first detach $I_1$ and $I_2$ from $T_1$.
Then $T_1\setminus(I_1\amalg I_2)$ is not necessarily connected.
We may assume that each connected component of $T_1\setminus(I_1\amalg I_2)$ is homeomorphic to a surface with exactly two (new) interchangeable holes (boundary components).
That is, we may assume that by detaching each $I_i$ from the surface, we have only one new hole in each component of the remaining surface, because in the case where there are two or more new holes, we can reduce them to one by increasing the genus (see Figure~\ref{fig:supp}).

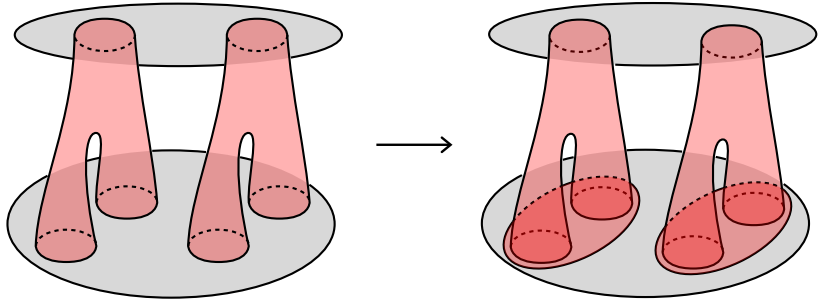
\begin{figure}[h]
    \centering
    \definecolor{cff0000}{RGB}{255,0,0}
    \def \globalscale {1.000000}
    \begin{tikzpicture}[y=0.80pt, x=0.80pt, yscale=-\globalscale, xscale=\globalscale, inner sep=0pt, outer sep=0pt]
      \path[fill=black,fill opacity=0.153] (66.7711,318.2960) .. controls (66.7711,299.0470) and (101.4000,283.4430) .. (144.1160,283.4430) .. controls (186.8330,283.4430) and (221.4620,299.0470) .. (221.4620,318.2960) .. controls (221.4620,337.5450) and (186.8330,353.1490) .. (144.1160,353.1490) .. controls (101.4000,353.1490) and (66.7711,337.5450) .. (66.7711,318.2960) -- cycle;

      \path[fill=black,fill opacity=0.151] (70.2175,228.8040) .. controls (70.2175,220.6390) and (104.8460,214.0200) .. (147.5630,214.0200) .. controls (190.2800,214.0200) and (224.9080,220.6390) .. (224.9080,228.8040) .. controls (224.9080,236.9690) and (190.2800,243.5880) .. (147.5630,243.5880) .. controls (104.8460,243.5880) and (70.2175,236.9690) .. (70.2175,228.8040) -- cycle;

      \path[draw=black,line width=0.800pt] (209.5210,299.6830) .. controls (217.0840,305.0680) and (221.4620,311.4520) .. (221.4620,318.2960) .. controls (221.4620,337.5450) and (186.8330,353.1490) .. (144.1160,353.1490) .. controls (101.4000,353.1490) and (66.7711,337.5450) .. (66.7711,318.2960) .. controls (66.7711,309.5580) and (73.9073,301.5710) .. (85.6968,295.4540);

      \path[draw=black,line width=0.800pt] (176.7400,286.6860) .. controls (178.5990,287.0760) and (180.4180,287.4990) .. (182.1930,287.9520);

      \path[draw=black,line width=0.800pt] (135.2200,283.6710) .. controls (138.1390,283.5200) and (141.1080,283.4430) .. (144.1160,283.4430) .. controls (149.7910,283.4430) and (155.3220,283.7180) .. (160.6490,284.2410);

      \path[draw=black,line width=0.800pt] (104.3240,288.4030) .. controls (106.3350,287.8580) and (108.4060,287.3540) .. (110.5300,286.8910);

      \path[draw=black,line width=0.800pt] (169.3740,242.9920) .. controls (162.4570,243.3790) and (155.1370,243.5880) .. (147.5630,243.5880) .. controls (141.0810,243.5880) and (134.7850,243.4350) .. (128.7690,243.1480);

      \path[draw=black,line width=0.800pt] (97.7578,240.1150) .. controls (80.9170,237.4030) and (70.2175,233.3420) .. (70.2175,228.8040) .. controls (70.2175,220.6390) and (104.8460,214.0200) .. (147.5630,214.0200) .. controls (190.2800,214.0200) and (224.9080,220.6390) .. (224.9080,228.8040) .. controls (224.9080,233.0280) and (215.6410,236.8380) .. (200.7810,239.5320);

      \path[fill=black,fill opacity=0.153] (291.0660,317.9980) .. controls (291.0660,298.7490) and (325.6950,283.1450) .. (368.4120,283.1450) .. controls (411.1280,283.1450) and (445.7570,298.7490) .. (445.7570,317.9980) .. controls (445.7570,337.2470) and (411.1280,352.8510) .. (368.4120,352.8510) .. controls (325.6950,352.8510) and (291.0660,337.2470) .. (291.0660,317.9980) -- cycle;

      \path[fill=black,fill opacity=0.151] (294.5130,228.5060) .. controls (294.5130,220.3410) and (329.1410,213.7220) .. (371.8580,213.7220) .. controls (414.5750,213.7220) and (449.2030,220.3410) .. (449.2030,228.5060) .. controls (449.2030,236.6700) and (414.5750,243.2890) .. (371.8580,243.2890) .. controls (329.1410,243.2890) and (294.5130,236.6700) .. (294.5130,228.5060) -- cycle;

      \path[draw=black,line width=0.800pt] (433.8160,299.3850) .. controls (441.3790,304.7700) and (445.7570,311.1540) .. (445.7570,317.9980) .. controls (445.7570,337.2470) and (411.1280,352.8510) .. (368.4120,352.8510) .. controls (325.6950,352.8510) and (291.0660,337.2470) .. (291.0660,317.9980) .. controls (291.0660,309.2600) and (298.2020,301.2730) .. (309.9920,295.1550);

      \path[draw=black,line width=0.800pt] (322.0530,239.8170) .. controls (305.2120,237.1050) and (294.5130,233.0440) .. (294.5130,228.5060) .. controls (294.5130,220.3410) and (329.1410,213.7220) .. (371.8580,213.7220) .. controls (414.5750,213.7220) and (449.2030,220.3410) .. (449.2030,228.5060) .. controls (449.2030,232.7290) and (439.9360,236.5390) .. (425.0770,239.2340);

      \path[draw=black,line width=0.800pt] (401.6890,286.5280) .. controls (403.5490,286.9190) and (404.7130,287.2000) .. (406.4890,287.6540);

      \path[draw=black,line width=0.800pt] (359.5150,283.3730) .. controls (362.4340,283.2220) and (365.4030,283.1450) .. (368.4120,283.1450) .. controls (374.0860,283.1450) and (380.4260,283.5020) .. (385.7530,284.0250);

      \path[draw=black,line width=0.800pt] (329.3440,287.9140) .. controls (331.3550,287.3690) and (332.3430,287.1360) .. (334.4670,286.6740);

      \path[draw=black,line width=0.800pt] (393.6700,242.6930) .. controls (386.7520,243.0810) and (379.4320,243.2890) .. (371.8580,243.2890) .. controls (365.3760,243.2890) and (360.1620,243.1820) .. (354.1450,242.8950);

      \path[draw=black,fill=cff0000,fill opacity=0.304,line width=0.800pt] (170.5060,228.8040) .. controls (170.5060,268.4500) and (153.1590,301.7610) .. (152.1970,328.6270) .. controls (152.1970,328.6270) and (152.3370,336.2660) .. (166.6610,336.2660) .. controls (180.9840,336.2660) and (180.8440,328.6270) .. (180.8440,328.6270) .. controls (174.2730,300.4580) and (174.3850,274.6590) .. (180.8440,275.2080) .. controls (186.0320,275.6480) and (180.2650,296.0440) .. (180.8440,307.8010) .. controls (180.8440,307.8010) and (180.9840,315.7360) .. (195.3070,315.7360) .. controls (209.6300,315.7360) and (209.6300,308.0970) .. (209.6300,308.0970) .. controls (207.7680,284.6180) and (201.0770,257.4570) .. (199.0340,228.8040) .. controls (199.0340,228.8040) and (199.0340,221.1650) .. (184.8290,221.1650) .. controls (170.5060,221.1650) and (170.5060,228.8040) .. (170.5060,228.8040) -- cycle;

      \path[draw=black,dash pattern=on 1.60pt,line width=0.800pt] (152.3370,328.6270) .. controls (152.3370,324.4080) and (158.7500,320.9880) .. (166.6610,320.9880) .. controls (174.5710,320.9880) and (180.9840,324.4080) .. (180.9840,328.6270);

      \path[draw=black,dash pattern=on 1.60pt,line width=0.800pt] (180.9840,308.0970) .. controls (180.9840,303.8780) and (187.3960,300.4580) .. (195.3070,300.4580) .. controls (203.2170,300.4580) and (209.6300,303.8780) .. (209.6300,308.0970);

      \path[draw=black,dash pattern=on 1.60pt,line width=0.800pt] (199.1520,228.8040) .. controls (199.1520,233.0230) and (192.7400,236.4430) .. (184.8290,236.4430) .. controls (176.9190,236.4430) and (170.5060,233.0230) .. (170.5060,228.8040);

      \path[draw=black,fill=cff0000,fill opacity=0.301,line width=0.800pt] (98.4414,228.8040) .. controls (98.4414,268.4500) and (81.0942,301.7610) .. (80.1328,328.6270) .. controls (80.1328,328.6270) and (80.2727,336.2660) .. (94.5959,336.2660) .. controls (108.9190,336.2660) and (108.7790,328.6270) .. (108.7790,328.6270) .. controls (102.2090,300.4580) and (102.3210,274.6590) .. (108.7790,275.2080) .. controls (113.9670,275.6480) and (108.2000,296.0440) .. (108.7790,307.8010) .. controls (108.7790,307.8010) and (108.9190,315.7360) .. (123.2420,315.7360) .. controls (137.5660,315.7360) and (137.5660,308.0970) .. (137.5660,308.0970) .. controls (135.7040,284.6180) and (129.0120,257.4570) .. (126.9690,228.8040) .. controls (126.9690,228.8040) and (126.9690,221.1650) .. (112.7650,221.1650) .. controls (98.4414,221.1650) and (98.4414,228.8040) .. (98.4414,228.8040) -- cycle;

      \path[draw=black,dash pattern=on 1.60pt,line width=0.800pt] (80.2727,328.6270) .. controls (80.2727,324.4080) and (86.6854,320.9880) .. (94.5959,320.9880) .. controls (102.5060,320.9880) and (108.9190,324.4080) .. (108.9190,328.6270);

      \path[draw=black,dash pattern=on 1.60pt,line width=0.800pt] (108.9190,308.0970) .. controls (108.9190,303.8780) and (115.3320,300.4580) .. (123.2420,300.4580) .. controls (131.1530,300.4580) and (137.5660,303.8780) .. (137.5660,308.0970);

      \path[draw=black,dash pattern=on 1.60pt,line width=0.800pt] (127.0880,228.8040) .. controls (127.0880,233.0230) and (120.6750,236.4430) .. (112.7650,236.4430) .. controls (104.8540,236.4430) and (98.4414,233.0230) .. (98.4414,228.8040);

      \path[draw=black,dash pattern=on 1.60pt,line width=0.800pt] (376.6320,331.7180) .. controls (376.6320,327.4990) and (383.0450,324.0780) .. (390.9560,324.0780) .. controls (398.8660,324.0780) and (405.2790,327.4990) .. (405.2790,331.7180);

      \path[draw=black,dash pattern=on 1.60pt,line width=0.800pt] (405.2790,311.1880) .. controls (405.2790,306.9690) and (411.6920,303.5490) .. (419.6020,303.5490) .. controls (427.5130,303.5490) and (433.9250,306.9690) .. (433.9250,311.1880);

      \path[draw=black,dash pattern=on 1.60pt,line width=0.800pt] (423.4480,231.8950) .. controls (423.4480,236.1140) and (417.0350,239.5340) .. (409.1240,239.5340) .. controls (401.2140,239.5340) and (394.8010,236.1140) .. (394.8010,231.8950);

      \path[fill=cff0000,fill opacity=0.301] (394.8710,231.8950) .. controls (394.8710,271.5410) and (377.5240,304.8520) .. (376.5630,331.7180) .. controls (376.5630,331.7180) and (376.7020,339.3570) .. (391.0260,339.3570) .. controls (405.3490,339.3570) and (405.2090,331.7180) .. (405.2090,331.7180) .. controls (398.6380,303.5490) and (398.7510,277.7500) .. (405.2090,278.2980) .. controls (410.3970,278.7390) and (404.6300,299.1350) .. (405.2090,310.8920) .. controls (405.2090,310.8920) and (405.3490,318.8270) .. (419.6720,318.8270) .. controls (433.9950,318.8270) and (433.9950,311.1880) .. (433.9950,311.1880) .. controls (432.1330,287.7090) and (425.4420,260.5480) .. (423.3990,231.8950) .. controls (423.3990,231.8950) and (423.3990,224.2560) .. (409.1940,224.2560) .. controls (394.8710,224.2560) and (394.8710,231.8950) .. (394.8710,231.8950) -- cycle;

      \path[draw=black,line width=0.800pt] (394.8710,231.8950) .. controls (394.8710,271.5410) and (377.5240,304.8520) .. (376.5630,331.7180) .. controls (376.5630,331.7180) and (376.7020,339.3570) .. (391.0260,339.3570) .. controls (405.3490,339.3570) and (405.2090,331.7180) .. (405.2090,331.7180) .. controls (398.6380,303.5490) and (398.7510,277.7500) .. (405.2090,278.2980) .. controls (410.3970,278.7390) and (404.6300,299.1350) .. (405.2090,310.8920) .. controls (405.2090,310.8920) and (405.3490,318.8270) .. (419.6720,318.8270) .. controls (433.9950,318.8270) and (433.9950,311.1880) .. (433.9950,311.1880) .. controls (432.1330,287.7090) and (425.4420,260.5480) .. (423.3990,231.8950) .. controls (423.3990,231.8950) and (423.3990,224.2560) .. (409.1940,224.2560) .. controls (394.8710,224.2560) and (394.8710,231.8950) .. (394.8710,231.8950) -- cycle;

      \path[draw=black,line width=0.800pt] (433.4750,301.9360) .. controls (434.7560,302.9180) and (435.7560,304.1130) .. (436.4230,305.5170) .. controls (440.6030,314.3300) and (430.0480,328.0880) .. (412.8470,336.2460) .. controls (395.6470,344.4040) and (378.3150,343.8740) .. (374.1350,335.0610) .. controls (372.1070,330.7850) and (373.5480,325.3440) .. (377.5340,319.9440);

      \path[draw=black,dash pattern=on 1.60pt,line width=0.800pt] (405.1410,301.3360) .. controls (416.9450,297.4020) and (427.8830,297.6470) .. (433.4750,301.9360);

      \path[draw=black,line width=0.800pt] (401.0750,302.8460) .. controls (402.4380,302.2880) and (403.7950,301.7850) .. (405.1410,301.3360);

      \path[draw=black,dash pattern=on 1.60pt,line width=0.800pt] (377.5340,319.9440) .. controls (381.7620,314.2150) and (388.8560,308.5320) .. (397.7100,304.3320) .. controls (398.8330,303.8000) and (399.9560,303.3040) .. (401.0750,302.8460);

      \path[fill=cff0000,fill opacity=0.301] (374.1350,335.0610) .. controls (369.9550,326.2480) and (380.5100,312.4900) .. (397.7100,304.3320) .. controls (414.9110,296.1740) and (432.2430,296.7050) .. (436.4230,305.5170) .. controls (440.6030,314.3300) and (430.0480,328.0880) .. (412.8470,336.2460) .. controls (395.6470,344.4040) and (378.3150,343.8740) .. (374.1350,335.0610) -- cycle;

      \path[draw=black,dash pattern=on 1.60pt,line width=0.800pt] (304.7960,329.0400) .. controls (304.7960,324.8210) and (311.2090,321.4010) .. (319.1190,321.4010) .. controls (327.0300,321.4010) and (333.4430,324.8210) .. (333.4430,329.0400);

      \path[draw=black,dash pattern=on 1.60pt,line width=0.800pt] (333.4430,308.5100) .. controls (333.4430,304.2910) and (339.8550,300.8710) .. (347.7660,300.8710) .. controls (355.6760,300.8710) and (362.0890,304.2910) .. (362.0890,308.5100);

      \path[draw=black,dash pattern=on 1.60pt,line width=0.800pt] (351.6110,229.2170) .. controls (351.6110,233.4360) and (345.1990,236.8560) .. (337.2880,236.8560) .. controls (329.3780,236.8560) and (322.9650,233.4360) .. (322.9650,229.2170);

      \path[fill=cff0000,fill opacity=0.301] (323.0350,229.2170) .. controls (323.0350,268.8640) and (305.6880,302.1750) .. (304.7260,329.0400) .. controls (304.7260,329.0400) and (304.8660,336.6790) .. (319.1890,336.6790) .. controls (333.5130,336.6790) and (333.3730,329.0400) .. (333.3730,329.0400) .. controls (326.8020,300.8710) and (326.9140,275.0730) .. (333.3730,275.6210) .. controls (338.5610,276.0620) and (332.7930,296.4580) .. (333.3730,308.2140) .. controls (333.3730,308.2140) and (333.5130,316.1490) .. (347.8360,316.1490) .. controls (362.1590,316.1490) and (362.1590,308.5100) .. (362.1590,308.5100) .. controls (360.2970,285.0310) and (353.6050,257.8700) .. (351.5630,229.2170) .. controls (351.5630,229.2170) and (351.5630,221.5780) .. (337.3580,221.5780) .. controls (323.0350,221.5780) and (323.0350,229.2170) .. (323.0350,229.2170) -- cycle;

      \path[draw=black,line width=0.800pt] (323.0350,229.2170) .. controls (323.0350,268.8640) and (305.6880,302.1750) .. (304.7260,329.0400) .. controls (304.7260,329.0400) and (304.8660,336.6790) .. (319.1890,336.6790) .. controls (333.5130,336.6790) and (333.3730,329.0400) .. (333.3730,329.0400) .. controls (326.8020,300.8710) and (326.9140,275.0730) .. (333.3730,275.6210) .. controls (338.5610,276.0620) and (332.7930,296.4580) .. (333.3730,308.2140) .. controls (333.3730,308.2140) and (333.5130,316.1490) .. (347.8360,316.1490) .. controls (362.1590,316.1490) and (362.1590,308.5100) .. (362.1590,308.5100) .. controls (360.2970,285.0310) and (353.6050,257.8700) .. (351.5630,229.2170) .. controls (351.5630,229.2170) and (351.5630,221.5780) .. (337.3580,221.5780) .. controls (323.0350,221.5780) and (323.0350,229.2170) .. (323.0350,229.2170) -- cycle;

      \path[draw=black,line width=0.800pt] (361.6390,299.2590) .. controls (362.9200,300.2410) and (363.9200,301.4350) .. (364.5860,302.8400) .. controls (368.7660,311.6530) and (358.2110,325.4110) .. (341.0110,333.5690) .. controls (323.8110,341.7270) and (306.4790,341.1960) .. (302.2990,332.3830) .. controls (300.2710,328.1070) and (301.7120,322.6670) .. (305.6970,317.2670);

      \path[draw=black,dash pattern=on 1.60pt,line width=0.800pt] (333.3050,298.6590) .. controls (345.1090,294.7240) and (356.0460,294.9700) .. (361.6390,299.2590);

      \path[draw=black,line width=0.800pt] (329.2390,300.1690) .. controls (330.6010,299.6110) and (331.9590,299.1080) .. (333.3050,298.6590);

      \path[draw=black,dash pattern=on 1.60pt,line width=0.800pt] (305.6970,317.2670) .. controls (309.9260,311.5380) and (317.0200,305.8540) .. (325.8740,301.6550) .. controls (326.9960,301.1220) and (328.1190,300.6270) .. (329.2390,300.1690);

      \path[fill=cff0000,fill opacity=0.301] (302.2990,332.3830) .. controls (298.1190,323.5710) and (308.6740,309.8130) .. (325.8740,301.6550) .. controls (343.0740,293.4970) and (360.4060,294.0270) .. (364.5860,302.8400) .. controls (368.7660,311.6530) and (358.2110,325.4110) .. (341.0110,333.5690) .. controls (323.8110,341.7270) and (306.4790,341.1960) .. (302.2990,332.3830) -- cycle;

          \path[draw=black,line width=0.905pt] (275.7390,280.8000) .. controls (272.8350,280.8000) and (241.2000,280.8000) .. (241.2000,280.8000);

          \path[fill=black] (272.5140,284.6330) -- (276.6640,281.2380) -- (277.1990,280.8000) -- (276.6640,280.3620) -- (272.5140,276.9670) .. controls (272.2730,276.7700) and (271.9160,276.8050) .. (271.7180,277.0470) .. controls (271.5200,277.2890) and (271.5560,277.6450) .. (271.7980,277.8430) -- (275.9470,281.2380) -- (275.9470,280.3620) -- (271.7980,283.7570) .. controls (271.5560,283.9550) and (271.5200,284.3110) .. (271.7180,284.5530) .. controls (271.9160,284.7950) and (272.2730,284.8310) .. (272.5140,284.6330) -- cycle;

    \end{tikzpicture}
    \caption{Reducing two boundary components of $I_i$ (glued to a disk) to one.}
    \label{fig:supp}
\end{figure}

Although $T_1\setminus(I_1\amalg I_2)$ is not necessarily connected, we first consider the case where it is connected.
The disconnected case will be considered later.

Let two holes of $T_1\setminus(I_1\amalg I_2)$, generated by detaching $I_1,I_2$, collapse to two points.
Denote this new atomic surface by $\widetilde T_1$.
Then the embedding $\phi:B_n\hookrightarrow\Gamma_{g,k}$ collapses to another regular embedding $\phi':B_n\hookrightarrow\Gamma_{h,k}^{(n)}$, where $S_{h,k}=\cup_{i=1}^{n-1}\widetilde T_i$.

Now we consider $\phi':B_n\hookrightarrow\Gamma_{h,k}^{(n)},\beta_i\mapsto[\widetilde\tau_i]$,
where each $\widetilde\tau_i$ is a self-homeomorphism of $\widetilde T_i$ interchanging two marked points by a simple twist.
We claim that this embedding is equivalent to that induced from $d$-fold ($d\geq 2$) covering map over a disk with $n$ marked points.

Denote the atomic surface $\widetilde T_1$ by $S_{h,b}^{(2)}$ $(b\geq 1)$, that is, $\tau_1$ induces a simple twist $\widetilde\tau_1:S_{h,b}^{(2)}\rightarrow S_{h,b}^{(2)}$ interchanging two marked points $p_1$ and $p_2$, and fixing the boundary components pointwise.
The mapping class group $\Gamma_{h,b}$ may be regarded as a subgroup of the automorphism group of $\pi_1(S_{h,b})\cong F_{2h+b-1}$, that is, the isotopy class of a self-homeomorphism of $S_{h,b}$ is completely determined by its action on the fundamental group of the surface. 

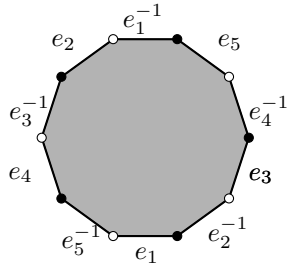
\begin{figure}[h]
    \centering
    \definecolor{cffffff}{RGB}{255,255,255}
    \def \globalscale {1.000000}
    \begin{tikzpicture}[y=0.80pt, x=0.80pt, yscale=-\globalscale, xscale=\globalscale, inner sep=0pt, outer sep=0pt]
      \path[draw=black,fill=black,fill opacity=0.301,line width=0.400pt] (145.4350,144.2240) -- (120.9340,161.9820) -- (90.6745,161.9480) -- (66.2141,144.1340) -- (56.8959,115.3450) -- (66.2792,86.5772) -- (90.7798,68.8188) -- (121.0390,68.8530) -- (145.5000,86.6669) -- (154.8180,115.4560) -- (145.4350,144.2240) -- cycle;

      \path[draw=black,line width=0.800pt] (66.2347,86.6240) -- (90.7259,68.8480);

      \path[draw=black,line width=0.800pt] (56.8694,115.4010) -- (66.2347,86.6240);

      \path[draw=black,line width=0.800pt] (66.2072,144.1860) -- (56.8694,115.4010);

      \path[draw=black,line width=0.800pt] (90.6814,161.9860) -- (66.2072,144.1860);

      \path[draw=black,line width=0.800pt] (120.9440,162.0000) -- (90.6814,161.9860);

      \path[draw=black,line width=0.800pt] (145.4350,144.2240) -- (120.9440,162.0000);

      \path[draw=black,line width=0.800pt] (154.8000,115.4470) -- (145.4350,144.2240);

      \path[draw=black,line width=0.800pt] (145.4620,86.6619) -- (154.8000,115.4470);

      \path[draw=black,line width=0.800pt] (120.9880,68.8625) -- (145.4620,86.6619);

      \path[draw=black,line width=0.800pt] (90.7259,68.8480) -- (120.9880,68.8625);

      \path[cm={{1.0,0.0,0.0,1.0,(100.944,167.514)}},fill=black] (0.0000,0.0000) node[below right] () {$e_1$};

      \path[cm={{1.0,0.0,0.0,1.0,(135.462,153.126)}},fill=black] (0.0000,0.0000) node[below right] () {$e_2^{-1}$};

      \path[cm={{1.0,0.0,0.0,1.0,(154.8,129.85)}},fill=black] (0.0000,0.0000) node[below right] () {$e_3$};

      \path[cm={{1.0,0.0,0.0,1.0,(154.8,95.5123)}},fill=black] (0.0000,0.0000) node[below right] () {$e_4^{-1}$};

      \path[cm={{1.0,0.0,0.0,1.0,(154.8,129.85)}},fill=black] (0.0000,0.0000) node[below right] () {$e_3$};

      \path[cm={{1.0,0.0,0.0,1.0,(140.117,66.736)}},fill=black] (0.0000,0.0000) node[below right] () {$e_5$};

      \path[cm={{1.0,0.0,0.0,1.0,(95.2692,52.1)}},fill=black] (0.0000,0.0000) node[below right] () {$e_1^{-1}$};

      \path[cm={{1.0,0.0,0.0,1.0,(66.2072,156.486)}},fill=black] (0.0000,0.0000) node[below right] () {$e_5^{-1}$};

      \path[cm={{1.0,0.0,0.0,1.0,(41.2,129.4)}},fill=black] (0.0000,0.0000) node[below right] () {$e_4$};

      \path[cm={{1.0,0.0,0.0,1.0,(41.5383,95.5123)}},fill=black] (0.0000,0.0000) node[below right] () {$e_3^{-1}$};

      \path[cm={{1.0,0.0,0.0,1.0,(61.5383,66.736)}},fill=black] (0.0000,0.0000) node[below right] () {$e_2$};

      \path[draw=black,fill=cffffff,line width=0.400pt] (88.5422,161.9820) .. controls (88.5422,160.8010) and (89.4999,159.8430) .. (90.6814,159.8430) .. controls (91.8628,159.8430) and (92.8205,160.8010) .. (92.8205,161.9820) .. controls (92.8205,163.1640) and (91.8628,164.1220) .. (90.6814,164.1220) .. controls (89.4999,164.1220) and (88.5422,163.1640) .. (88.5422,161.9820) -- cycle;

      \path[draw=black,fill=black,line width=0.400pt] (118.8490,161.9820) .. controls (118.8490,160.8010) and (119.8070,159.8430) .. (120.9880,159.8430) .. controls (122.1690,159.8430) and (123.1270,160.8010) .. (123.1270,161.9820) .. controls (123.1270,163.1640) and (122.1690,164.1220) .. (120.9880,164.1220) .. controls (119.8070,164.1220) and (118.8490,163.1640) .. (118.8490,161.9820) -- cycle;

      \path[draw=black,fill=cffffff,line width=0.400pt] (143.3230,144.1860) .. controls (143.3230,143.0050) and (144.2810,142.0470) .. (145.4620,142.0470) .. controls (146.6440,142.0470) and (147.6010,143.0050) .. (147.6010,144.1860) .. controls (147.6010,145.3680) and (146.6440,146.3250) .. (145.4620,146.3250) .. controls (144.2810,146.3250) and (143.3230,145.3680) .. (143.3230,144.1860) -- cycle;

      \path[draw=black,fill=black,line width=0.400pt] (152.6790,115.4620) .. controls (152.6790,114.2800) and (153.6370,113.3230) .. (154.8180,113.3230) .. controls (155.9990,113.3230) and (156.9570,114.2800) .. (156.9570,115.4620) .. controls (156.9570,116.6430) and (155.9990,117.6010) .. (154.8180,117.6010) .. controls (153.6370,117.6010) and (152.6790,116.6430) .. (152.6790,115.4620) -- cycle;

      \path[draw=black,fill=cffffff,line width=0.400pt] (143.3230,86.6240) .. controls (143.3230,85.4426) and (144.2810,84.4848) .. (145.4620,84.4848) .. controls (146.6440,84.4848) and (147.6010,85.4426) .. (147.6010,86.6240) .. controls (147.6010,87.8055) and (146.6440,88.7632) .. (145.4620,88.7632) .. controls (144.2810,88.7632) and (143.3230,87.8055) .. (143.3230,86.6240) -- cycle;

      \path[draw=black,fill=black,line width=0.400pt] (118.8490,68.8188) .. controls (118.8490,67.6374) and (119.8070,66.6796) .. (120.9880,66.6796) .. controls (122.1690,66.6796) and (123.1270,67.6374) .. (123.1270,68.8188) .. controls (123.1270,70.0002) and (122.1690,70.9579) .. (120.9880,70.9579) .. controls (119.8070,70.9579) and (118.8490,70.0002) .. (118.8490,68.8188) -- cycle;

      \path[draw=black,fill=cffffff,line width=0.400pt] (88.5422,68.8188) .. controls (88.5422,67.6374) and (89.4999,66.6796) .. (90.6814,66.6796) .. controls (91.8628,66.6796) and (92.8205,67.6374) .. (92.8205,68.8188) .. controls (92.8205,70.0002) and (91.8628,70.9579) .. (90.6814,70.9579) .. controls (89.4999,70.9579) and (88.5422,70.0002) .. (88.5422,68.8188) -- cycle;

      \path[draw=black,fill=black,line width=0.400pt] (64.0680,86.6240) .. controls (64.0680,85.4426) and (65.0258,84.4848) .. (66.2072,84.4848) .. controls (67.3886,84.4848) and (68.3464,85.4426) .. (68.3464,86.6240) .. controls (68.3464,87.8055) and (67.3886,88.7632) .. (66.2072,88.7632) .. controls (65.0258,88.7632) and (64.0680,87.8055) .. (64.0680,86.6240) -- cycle;

      \path[draw=black,fill=cffffff,line width=0.400pt] (54.7567,115.4620) .. controls (54.7567,114.2800) and (55.7145,113.3230) .. (56.8959,113.3230) .. controls (58.0773,113.3230) and (59.0351,114.2800) .. (59.0351,115.4620) .. controls (59.0351,116.6430) and (58.0773,117.6010) .. (56.8959,117.6010) .. controls (55.7145,117.6010) and (54.7567,116.6430) .. (54.7567,115.4620) -- cycle;

      \path[draw=black,fill=black,line width=0.400pt] (64.0680,144.2380) .. controls (64.0680,143.0570) and (65.0258,142.0990) .. (66.2072,142.0990) .. controls (67.3886,142.0990) and (68.3464,143.0570) .. (68.3464,144.2380) .. controls (68.3464,145.4200) and (67.3886,146.3780) .. (66.2072,146.3780) .. controls (65.0258,146.3780) and (64.0680,145.4200) .. (64.0680,144.2380) -- cycle;

    \end{tikzpicture}
    \caption{$S_2$ as a 10-gon. The white dots represent $p_1$ and the black dots represent $p_2$.}
    \label{fig:polygon}
\end{figure}
On the other hand, it is well-known that closed surface $S_h$ can be expressed as a $(4h+2)$-gon with two vertices and opposite sides identified (Figure~\ref{fig:polygon}).
We may regard $S_{h,1}= S_h - D$ as the groupoid $\G$ with two vertices $p_1,p_2$ and $2h+1$ edges $e_1,e_2,\ldots,e_{2h+1}$ from $p_1$ to $p_2$.
$S_{h,b}$ has additional $b-1$ edges $e_{2h+2},\ldots,e_{2h+b}$ and a self-homeomorphism of $S_{h,b}$ is a self-functor fixing the loops along the boundary components.

\begin{figure}[h]
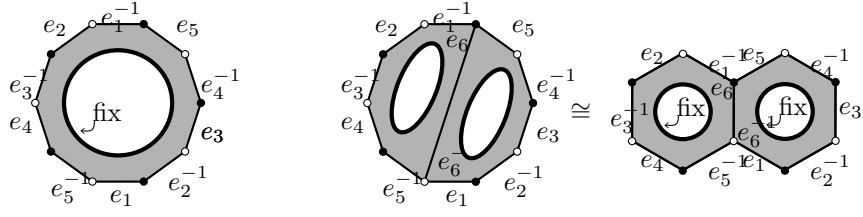

    \centering
    \definecolor{cffffff}{RGB}{255,255,255}
    \def \globalscale {0.8}

    \caption{$S_{2,1}$ as a 10-gon, and $S_{2,2}$ as two 6-gons.}
    \label{fig:polygon_b}
\end{figure}

Let us now figure out how $\widetilde\tau_1$ acts on $\G\simeq S_{h,b}^{(2)}$ as a self-functor.

\begin{enumerate}
    \item[Case I.] 
    In the case of $b=1$, a loop along the boundary (in counter clockwise direction) of the surface can be expressed as
    $$\partial = e_1\cdot e_2^{-1}\cdot\cdots\cdot e_{2h}^{-1}\cdot e_{2h+1}\cdot e_1^{-1}\cdot e_2\cdot\cdots\cdot e_{2h+1}^{-1}.$$
    $\widetilde\tau_1$ should move every edge of $\G$ while fixing the boundary loop $\partial$.
    By the definition of simple twist, each edge should be moved to its adjacent edge, that is, if we have the twist in the counter clockwise direction in $\G$ then $\widetilde\tau_1$ acts as:
    $$e_1\mapsto e_2^{-1},\ e_2^{-1}\mapsto e_3,\ldots,\ e_{2h}^{-1}\mapsto e_{2h+1},\ e_{2h+1}\mapsto e_1^{-1},$$
    $$e_1^{-1}\mapsto e_2,\ e_2\mapsto e_3^{-1},\ldots,\ e_{2h+1}^{-1}\mapsto e_1$$
    because $p_1$ and $p_2$ should be interchanged only once (see Figure~\ref{fig:tau_1}).
    \ \\

    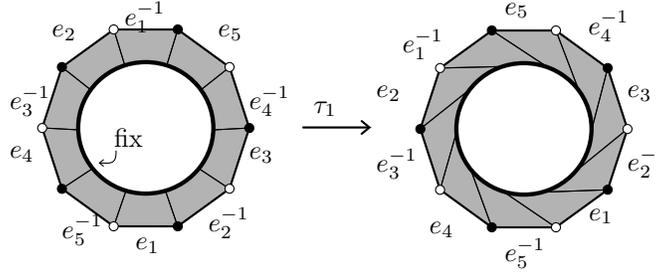
\begin{figure}[h]
        \centering
        \definecolor{cffffff}{RGB}{255,255,255}
        \def \globalscale {1.000000}
        \begin{tikzpicture}[y=0.80pt, x=0.80pt, yscale=-\globalscale, xscale=\globalscale, inner sep=0pt, outer sep=0pt]
          \path[draw=black,fill=black,fill opacity=0.301,line width=0.400pt] (324.2340,144.7000) -- (299.7330,162.4580) -- (269.4740,162.4240) -- (245.0130,144.6100) -- (235.6950,115.8210) -- (245.0780,87.0530) -- (269.5790,69.2946) -- (299.8390,69.3288) -- (324.2990,87.1427) -- (333.6170,115.9320) -- (324.2340,144.7000) -- cycle;

            \path[draw=black,fill=black,line width=0.400pt] (90.7259,68.8480) -- (95.9606,85.8794);

            \path[draw=black,fill=black,line width=0.400pt] (120.9880,68.8625) -- (115.8930,85.4329);

            \path[draw=black,fill=black,line width=0.400pt] (145.4620,86.6619) -- (131.9030,97.4167);

            \path[draw=black,fill=black,line width=0.400pt] (154.8000,115.4470) -- (137.7220,116.3650);

            \path[draw=black,fill=black,line width=0.400pt] (145.4350,144.2240) -- (131.1160,134.7650);

            \path[draw=black,fill=black,line width=0.400pt] (120.9440,162.0000) -- (115.2130,145.3160);

            \path[draw=black,fill=black,line width=0.400pt] (90.6814,161.9860) -- (96.4413,145.3200);

            \path[draw=black,fill=black,line width=0.400pt] (66.2072,144.1860) -- (80.6886,134.6500);

            \path[draw=black,fill=black,line width=0.400pt] (56.8694,115.4010) -- (74.1474,116.3450);

            \path[draw=black,fill=black,line width=0.400pt] (66.2347,86.6240) -- (79.6132,97.2441);

          \path[draw=black,fill=black,fill opacity=0.301,line width=0.400pt] (145.4350,144.2240) -- (120.9340,161.9820) -- (90.6745,161.9480) -- (66.2141,144.1340) -- (56.8959,115.3450) -- (66.2792,86.5772) -- (90.7798,68.8188) -- (121.0390,68.8530) -- (145.5000,86.6669) -- (154.8180,115.4560) -- (145.4350,144.2240) -- cycle;

          \path[draw=black,line width=0.800pt] (66.2347,86.6240) -- (90.7259,68.8480);

          \path[draw=black,line width=0.800pt] (56.8694,115.4010) -- (66.2347,86.6240);

          \path[draw=black,line width=0.800pt] (66.2072,144.1860) -- (56.8694,115.4010);

          \path[draw=black,line width=0.800pt] (90.6814,161.9860) -- (66.2072,144.1860);

          \path[draw=black,line width=0.800pt] (120.9440,162.0000) -- (90.6814,161.9860);

          \path[draw=black,line width=0.800pt] (145.4350,144.2240) -- (120.9440,162.0000);

          \path[draw=black,line width=0.800pt] (154.8000,115.4470) -- (145.4350,144.2240);

          \path[draw=black,line width=0.800pt] (145.4620,86.6619) -- (154.8000,115.4470);

          \path[draw=black,line width=0.800pt] (120.9880,68.8625) -- (145.4620,86.6619);

          \path[draw=black,line width=0.800pt] (90.7259,68.8480) -- (120.9880,68.8625);

          \path[cm={{1.0,0.0,0.0,1.0,(100.944,167.514)}},fill=black] (0.0000,0.0000) node[below right] () {$e_{1}$};

          \path[cm={{1.0,0.0,0.0,1.0,(135.462,153.126)}},fill=black] (0.0000,0.0000) node[below right] () {$e_{2}^{-1}$};

          \path[cm={{1.0,0.0,0.0,1.0,(154.8,95.5123)}},fill=black] (0.0000,0.0000) node[below right] () {$e_{4}^{-1}$};

          \path[cm={{1.0,0.0,0.0,1.0,(154.8,124.769)}},fill=black] (0.0000,0.0000) node[below right] () {$e_{3}$};

          \path[cm={{1.0,0.0,0.0,1.0,(140.117,66.736)}},fill=black] (0.0000,0.0000) node[below right] () {$e_{5}$};

          \path[cm={{1.0,0.0,0.0,1.0,(95.7588,55.736)}},fill=black] (0.0000,0.0000) node[below right] () {$e_{1}^{-1}$};

          \path[cm={{1.0,0.0,0.0,1.0,(66.2072,156.486)}},fill=black] (0.0000,0.0000) node[below right] () {$e_{5}^{-1}$};

          \path[cm={{1.0,0.0,0.0,1.0,(41.5383,124.35)}},fill=black] (0.0000,0.0000) node[below right] () {$e_{4}$};

          \path[cm={{1.0,0.0,0.0,1.0,(41.5383,95.5123)}},fill=black] (0.0000,0.0000) node[below right] () {$e_{3}^{-1}$};

          \path[cm={{1.0,0.0,0.0,1.0,(61.5383,66.736)}},fill=black] (0.0000,0.0000) node[below right] () {$e_{2}$};

          \path[draw=black,fill=cffffff,line width=1.600pt] (73.9019,115.4010) .. controls (73.9019,98.2231) and (88.2087,84.2980) .. (105.8570,84.2980) .. controls (123.5050,84.2980) and (137.8120,98.2231) .. (137.8120,115.4010) .. controls (137.8120,132.5780) and (123.5050,146.5030) .. (105.8570,146.5030) .. controls (88.2087,146.5030) and (73.9019,132.5780) .. (73.9019,115.4010) -- cycle;

          \path[cm={{1.0,0.0,0.0,1.0,(90.9436,115.462)}},fill=black] (0.0000,0.0000) node[below right] () {fix};

              \path[draw=black,line width=0.400pt] (84.2090,131.9840) .. controls (88.5738,132.0940) and (92.1916,133.1130) .. (91.7775,126.1410);

              \path[fill=black] (85.6565,130.3090) -- (83.8033,131.7850) -- (83.5643,131.9750) -- (83.7982,132.1720) -- (85.6114,133.6960) .. controls (85.7171,133.7850) and (85.8748,133.7710) .. (85.9636,133.6650) .. controls (86.0525,133.5600) and (86.0388,133.4020) .. (85.9331,133.3130) -- (84.1199,131.7890) -- (84.1148,132.1760) -- (85.9679,130.7000) .. controls (86.0759,130.6140) and (86.0938,130.4570) .. (86.0078,130.3490) .. controls (85.9218,130.2410) and (85.7645,130.2230) .. (85.6565,130.3090) -- cycle;

          \path[draw=black,fill=cffffff,line width=0.400pt] (88.5422,161.9820) .. controls (88.5422,160.8010) and (89.4999,159.8430) .. (90.6814,159.8430) .. controls (91.8628,159.8430) and (92.8205,160.8010) .. (92.8205,161.9820) .. controls (92.8205,163.1640) and (91.8628,164.1220) .. (90.6814,164.1220) .. controls (89.4999,164.1220) and (88.5422,163.1640) .. (88.5422,161.9820) -- cycle;

          \path[draw=black,fill=black,line width=0.400pt] (118.8490,161.9820) .. controls (118.8490,160.8010) and (119.8070,159.8430) .. (120.9880,159.8430) .. controls (122.1690,159.8430) and (123.1270,160.8010) .. (123.1270,161.9820) .. controls (123.1270,163.1640) and (122.1690,164.1220) .. (120.9880,164.1220) .. controls (119.8070,164.1220) and (118.8490,163.1640) .. (118.8490,161.9820) -- cycle;

          \path[draw=black,fill=cffffff,line width=0.400pt] (143.3230,144.1860) .. controls (143.3230,143.0050) and (144.2810,142.0470) .. (145.4620,142.0470) .. controls (146.6440,142.0470) and (147.6010,143.0050) .. (147.6010,144.1860) .. controls (147.6010,145.3680) and (146.6440,146.3250) .. (145.4620,146.3250) .. controls (144.2810,146.3250) and (143.3230,145.3680) .. (143.3230,144.1860) -- cycle;

          \path[draw=black,fill=black,line width=0.400pt] (152.6790,115.4620) .. controls (152.6790,114.2800) and (153.6370,113.3230) .. (154.8180,113.3230) .. controls (155.9990,113.3230) and (156.9570,114.2800) .. (156.9570,115.4620) .. controls (156.9570,116.6430) and (155.9990,117.6010) .. (154.8180,117.6010) .. controls (153.6370,117.6010) and (152.6790,116.6430) .. (152.6790,115.4620) -- cycle;

          \path[draw=black,fill=cffffff,line width=0.400pt] (143.3230,86.6240) .. controls (143.3230,85.4426) and (144.2810,84.4848) .. (145.4620,84.4848) .. controls (146.6440,84.4848) and (147.6010,85.4426) .. (147.6010,86.6240) .. controls (147.6010,87.8055) and (146.6440,88.7632) .. (145.4620,88.7632) .. controls (144.2810,88.7632) and (143.3230,87.8055) .. (143.3230,86.6240) -- cycle;

          \path[draw=black,fill=black,line width=0.400pt] (118.8490,68.8188) .. controls (118.8490,67.6374) and (119.8070,66.6796) .. (120.9880,66.6796) .. controls (122.1690,66.6796) and (123.1270,67.6374) .. (123.1270,68.8188) .. controls (123.1270,70.0002) and (122.1690,70.9579) .. (120.9880,70.9579) .. controls (119.8070,70.9579) and (118.8490,70.0002) .. (118.8490,68.8188) -- cycle;

          \path[draw=black,fill=cffffff,line width=0.400pt] (88.5422,68.8188) .. controls (88.5422,67.6374) and (89.4999,66.6796) .. (90.6814,66.6796) .. controls (91.8628,66.6796) and (92.8205,67.6374) .. (92.8205,68.8188) .. controls (92.8205,70.0002) and (91.8628,70.9579) .. (90.6814,70.9579) .. controls (89.4999,70.9579) and (88.5422,70.0002) .. (88.5422,68.8188) -- cycle;

          \path[draw=black,fill=black,line width=0.400pt] (64.0680,86.6240) .. controls (64.0680,85.4426) and (65.0258,84.4848) .. (66.2072,84.4848) .. controls (67.3886,84.4848) and (68.3464,85.4426) .. (68.3464,86.6240) .. controls (68.3464,87.8055) and (67.3886,88.7632) .. (66.2072,88.7632) .. controls (65.0258,88.7632) and (64.0680,87.8055) .. (64.0680,86.6240) -- cycle;

          \path[draw=black,fill=cffffff,line width=0.400pt] (54.7567,115.4620) .. controls (54.7567,114.2800) and (55.7145,113.3230) .. (56.8959,113.3230) .. controls (58.0773,113.3230) and (59.0351,114.2800) .. (59.0351,115.4620) .. controls (59.0351,116.6430) and (58.0773,117.6010) .. (56.8959,117.6010) .. controls (55.7145,117.6010) and (54.7567,116.6430) .. (54.7567,115.4620) -- cycle;

          \path[draw=black,fill=black,line width=0.400pt] (64.0680,144.2380) .. controls (64.0680,143.0570) and (65.0258,142.0990) .. (66.2072,142.0990) .. controls (67.3886,142.0990) and (68.3464,143.0570) .. (68.3464,144.2380) .. controls (68.3464,145.4200) and (67.3886,146.3780) .. (66.2072,146.3780) .. controls (65.0258,146.3780) and (64.0680,145.4200) .. (64.0680,144.2380) -- cycle;

          \path[draw=black,fill=black,line width=0.400pt] (245.0340,87.0998) -- (274.7600,86.3552);

          \path[draw=black,fill=black,line width=0.400pt] (269.5250,69.3238) -- (294.6930,85.9088);

          \path[draw=black,fill=black,line width=0.400pt] (299.7870,69.3383) -- (310.7030,97.8925);

          \path[draw=black,fill=black,line width=0.400pt] (324.2610,87.1377) -- (316.5210,116.8410);

          \path[draw=black,fill=black,line width=0.400pt] (333.5990,115.9230) -- (309.9150,135.2400);

          \path[draw=black,fill=black,line width=0.400pt] (324.2340,144.7000) -- (294.0120,145.7920);

          \path[draw=black,fill=black,line width=0.400pt] (299.7430,162.4760) -- (275.2400,145.7960);

          \path[draw=black,fill=black,line width=0.400pt] (269.4800,162.4610) -- (259.4880,135.1260);

          \path[draw=black,fill=black,line width=0.400pt] (245.0060,144.6620) -- (252.9470,116.8200);

          \path[draw=black,fill=black,line width=0.400pt] (235.6690,115.8760) -- (258.4120,97.7199);

          \path[draw=black,line width=0.800pt] (245.0340,87.0998) -- (269.5250,69.3238);

          \path[draw=black,line width=0.800pt] (235.6690,115.8760) -- (245.0340,87.0998);

          \path[draw=black,line width=0.800pt] (245.0060,144.6620) -- (235.6690,115.8760);

          \path[draw=black,line width=0.800pt] (269.4800,162.4610) -- (245.0060,144.6620);

          \path[draw=black,line width=0.800pt] (299.7430,162.4760) -- (269.4800,162.4610);

          \path[draw=black,line width=0.800pt] (324.2340,144.7000) -- (299.7430,162.4760);

          \path[draw=black,line width=0.800pt] (333.5990,115.9230) -- (324.2340,144.7000);

          \path[draw=black,line width=0.800pt] (324.2610,87.1377) -- (333.5990,115.9230);

          \path[draw=black,line width=0.800pt] (299.7870,69.3383) -- (324.2610,87.1377);

          \path[draw=black,line width=0.800pt] (269.5250,69.3238) -- (299.7870,69.3383);

          \path[draw=black,fill=cffffff,line width=1.600pt] (252.7010,115.8760) .. controls (252.7010,98.6989) and (267.0080,84.7738) .. (284.6560,84.7738) .. controls (302.3040,84.7738) and (316.6110,98.6989) .. (316.6110,115.8760) .. controls (316.6110,133.0540) and (302.3040,146.9790) .. (284.6560,146.9790) .. controls (267.0080,146.9790) and (252.7010,133.0540) .. (252.7010,115.8760) -- cycle;

          \path[draw=black,fill=black,line width=0.400pt] (267.3410,162.4580) .. controls (267.3410,161.2770) and (268.2990,160.3190) .. (269.4800,160.3190) .. controls (270.6620,160.3190) and (271.6200,161.2770) .. (271.6200,162.4580) .. controls (271.6200,163.6400) and (270.6620,164.5970) .. (269.4800,164.5970) .. controls (268.2990,164.5970) and (267.3410,163.6400) .. (267.3410,162.4580) -- cycle;

          \path[draw=black,fill=cffffff,line width=0.400pt] (297.6480,162.4580) .. controls (297.6480,161.2770) and (298.6060,160.3190) .. (299.7870,160.3190) .. controls (300.9690,160.3190) and (301.9260,161.2770) .. (301.9260,162.4580) .. controls (301.9260,163.6400) and (300.9690,164.5970) .. (299.7870,164.5970) .. controls (298.6060,164.5970) and (297.6480,163.6400) .. (297.6480,162.4580) -- cycle;

          \path[draw=black,fill=black,line width=0.400pt] (322.1220,144.6620) .. controls (322.1220,143.4810) and (323.0800,142.5230) .. (324.2610,142.5230) .. controls (325.4430,142.5230) and (326.4010,143.4810) .. (326.4010,144.6620) .. controls (326.4010,145.8430) and (325.4430,146.8010) .. (324.2610,146.8010) .. controls (323.0800,146.8010) and (322.1220,145.8430) .. (322.1220,144.6620) -- cycle;

          \path[draw=black,fill=cffffff,line width=0.400pt] (331.4780,115.9380) .. controls (331.4780,114.7560) and (332.4360,113.7990) .. (333.6170,113.7990) .. controls (334.7990,113.7990) and (335.7560,114.7560) .. (335.7560,115.9380) .. controls (335.7560,117.1190) and (334.7990,118.0770) .. (333.6170,118.0770) .. controls (332.4360,118.0770) and (331.4780,117.1190) .. (331.4780,115.9380) -- cycle;

          \path[draw=black,fill=black,line width=0.400pt] (322.1220,87.0998) .. controls (322.1220,85.9184) and (323.0800,84.9607) .. (324.2610,84.9607) .. controls (325.4430,84.9607) and (326.4010,85.9184) .. (326.4010,87.0998) .. controls (326.4010,88.2813) and (325.4430,89.2390) .. (324.2610,89.2390) .. controls (323.0800,89.2390) and (322.1220,88.2813) .. (322.1220,87.0998) -- cycle;

          \path[draw=black,fill=cffffff,line width=0.400pt] (297.6480,69.2946) .. controls (297.6480,68.1132) and (298.6060,67.1554) .. (299.7870,67.1554) .. controls (300.9690,67.1554) and (301.9260,68.1132) .. (301.9260,69.2946) .. controls (301.9260,70.4760) and (300.9690,71.4338) .. (299.7870,71.4338) .. controls (298.6060,71.4338) and (297.6480,70.4760) .. (297.6480,69.2946) -- cycle;

          \path[draw=black,fill=black,line width=0.400pt] (267.3410,69.2946) .. controls (267.3410,68.1132) and (268.2990,67.1554) .. (269.4800,67.1554) .. controls (270.6620,67.1554) and (271.6200,68.1132) .. (271.6200,69.2946) .. controls (271.6200,70.4760) and (270.6620,71.4338) .. (269.4800,71.4338) .. controls (268.2990,71.4338) and (267.3410,70.4760) .. (267.3410,69.2946) -- cycle;

          \path[draw=black,fill=cffffff,line width=0.400pt] (242.8670,87.0998) .. controls (242.8670,85.9184) and (243.8250,84.9607) .. (245.0060,84.9607) .. controls (246.1880,84.9607) and (247.1460,85.9184) .. (247.1460,87.0998) .. controls (247.1460,88.2813) and (246.1880,89.2390) .. (245.0060,89.2390) .. controls (243.8250,89.2390) and (242.8670,88.2813) .. (242.8670,87.0998) -- cycle;

          \path[draw=black,fill=black,line width=0.400pt] (233.5560,115.9380) .. controls (233.5560,114.7560) and (234.5140,113.7990) .. (235.6950,113.7990) .. controls (236.8760,113.7990) and (237.8340,114.7560) .. (237.8340,115.9380) .. controls (237.8340,117.1190) and (236.8760,118.0770) .. (235.6950,118.0770) .. controls (234.5140,118.0770) and (233.5560,117.1190) .. (233.5560,115.9380) -- cycle;

          \path[draw=black,fill=cffffff,line width=0.400pt] (242.8670,144.7140) .. controls (242.8670,143.5330) and (243.8250,142.5750) .. (245.0060,142.5750) .. controls (246.1880,142.5750) and (247.1460,143.5330) .. (247.1460,144.7140) .. controls (247.1460,145.8960) and (246.1880,146.8530) .. (245.0060,146.8530) .. controls (243.8250,146.8530) and (242.8670,145.8960) .. (242.8670,144.7140) -- cycle;

            \path[fill=black] (180.0000,115.2000) -- (212.4000,115.2000);

              \path[draw=black,line width=0.800pt] (211.1100,115.2000) .. controls (208.5100,115.2000) and (180.0000,115.2000) .. (180.0000,115.2000);

              \path[fill=black] (208.2600,118.5870) -- (211.9270,115.5870) -- (212.4000,115.2000) -- (211.9270,114.8130) -- (208.2600,111.8130) .. controls (208.0460,111.6380) and (207.7310,111.6700) .. (207.5570,111.8830) .. controls (207.3820,112.0970) and (207.4130,112.4120) .. (207.6270,112.5870) -- (211.2940,115.5870) -- (211.2940,114.8130) -- (207.6270,117.8130) .. controls (207.4130,117.9880) and (207.3820,118.3030) .. (207.5570,118.5170) .. controls (207.7310,118.7300) and (208.0460,118.7620) .. (208.2600,118.5870) -- cycle;

          \path[cm={{1.0,0.0,0.0,1.0,(185.139,102.799)}},fill=black] (0.0000,0.0000) node[below right] () {$\tau_1$};

          \path[cm={{1.0,0.0,0.0,1.0,(315.35,153.562)}},fill=black] (0.0000,0.0000) node[below right] () {$e_{1}$};

          \path[cm={{1.0,0.0,0.0,1.0,(333.617,124.769)}},fill=black] (0.0000,0.0000) node[below right] () {$e_{2}^{-1}$};

          \path[cm={{1.0,0.0,0.0,1.0,(333.617,95.5123)}},fill=black] (0.0000,0.0000) node[below right] () {$e_{3}$};

          \path[cm={{1.0,0.0,0.0,1.0,(315.0,60.0)}},fill=black] (0.0000,0.0000) node[below right] () {$e_{4}^{-1}$};

          \path[cm={{1.0,0.0,0.0,1.0,(275.483,56.1554)}},fill=black] (0.0000,0.0000) node[below right] () {$e_{5}$};

          \path[cm={{1.0,0.0,0.0,1.0,(227.241,68.8188)}},fill=black] (0.0000,0.0000) node[below right] () {$e_{1}^{-1}$};

          \path[cm={{1.0,0.0,0.0,1.0,(214.917,95.5123)}},fill=black] (0.0000,0.0000) node[below right] () {$e_{2}$};

          \path[cm={{1.0,0.0,0.0,1.0,(214.917,125.241)}},fill=black] (0.0000,0.0000) node[below right] () {$e_{3}^{-1}$};

          \path[cm={{1.0,0.0,0.0,1.0,(240.0,160.0)}},fill=black] (0.0000,0.0000) node[below right] () {$e_{4}$};

          \path[cm={{1.0,0.0,0.0,1.0,(275.483,167.486)}},fill=black] (0.0000,0.0000) node[below right] () {$e_{5}^{-1}$};

        \end{tikzpicture}
        \caption{The action of $\tau_1$ on $S_{2,1}$.}
        \label{fig:tau_1}
    \end{figure}
    
    \item[Case II.] In the case of $b=2$, the groupoid $\G$ has an additional edge $e_{2h+2}$ to the groupoid of Case I (see Figure~\ref{fig:polygon_b}).
    Then two boundary loops $\partial_1,\partial_2$ represented by
    $$\partial_1=e_1\cdot e_2^{-1}\cdot\cdots\cdot e_{2h}^{-1}\cdot e_{2h+1}\cdot e_{2h+2}^{-1},\quad
    \partial_2=e_1^{-1}\cdot e_2\cdot\cdots\cdot e_{2h}\cdot e_{2h+1}^{-1}\cdot e_{2h+2}$$
    are supposed to be fixed.
    We may regard this groupoid as two identical $(2h+2)$-gons glued together along the edge $e_{2h+2}$.
    As in Case I, the simple twist $\widetilde\tau_1$ moves each edge to its adjacent edge in each $(2h+2)$-gon.
    \item[Case III.] If $b\geq 3$, then we may add edges $e_{2h+2},e_{2h+3},\ldots$ to the groupoid $\G$ in Case I.
    Then moving each edge to its adjacent edge in one direction in each polygon causes a mismatch of sides, in particular, the continuity of $\widetilde\tau_1$ does not hold (see Figure~\ref{fig:mismatch}).
\end{enumerate}

\begin{figure}[h]
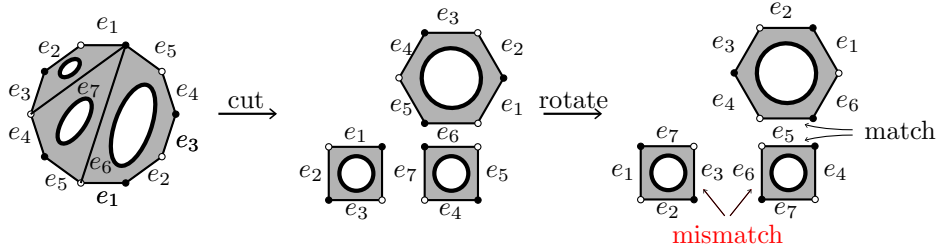

    \centering
    \definecolor{cffffff}{RGB}{255,255,255}
    \definecolor{cff0000}{RGB}{255,0,0}
    \def \globalscale {0.7}

    \caption{Mismatch of sides.}
    \label{fig:mismatch}
\end{figure}

In conclusion, the number of boundary components of the atomic surface $\widetilde T_1$ equals 1 or 2.
As a matter of fact, the twist $\widetilde\tau_1$ on $\widetilde T_1$ is equivalent to that induced by a $d$-fold ($d=2h+1$ or $2h+2$) cyclic branched covering over a disk with two marked points
because the simple twist $\widetilde\tau_1$ maps $e_i\mapsto e_{i+1}^{-1}$ where the indices read modulo $d$.
This can also be shown by the argument of ribbon graph automorphism (see Remark~\ref{rmk:ribbon}).

The embedding $\phi':B_n\hookrightarrow\Gamma_{g,k}^{(n)}$ is equivalent to that induced from $d$-fold covering over a disk with $n$ marked points.
Now we have
$$g=\frac{1}{2}(dn-n-d-\gcd(d,n))+1\quad\textrm{and}\quad k=\gcd(d,n)$$
by the Riemann-Hurwitz formula.
This means that if each $\widetilde T_i$ is connected, as we have seen in section~\ref{sec:d-fold},
we have a $\mathcal{D}$-algebra map $\Phi : \mathcal{C} \rightarrow \mathcal{S}$ induced by $\phi$, where
$\mathcal{C}_{m} = \conf_{dm}(D)$ and $\mathcal{S}_{m} = \mathcal{M}_{g(d,m), k}$.
Therefore, $\B\phi^{+} : \B B_{\infty}^{+} \rightarrow \B\Gamma_{\infty}^{+}$ is a double loop space map, which means the embedding $\phi:B_n\hookrightarrow\Gamma_{g,k}^{(n)}$ satisfies the homology triviality.

Finally, consider the case where $\widetilde T_1$ is not connected.
Let $\widetilde T_1 = S_1\amalg\cdots\amalg S_m$, where $S_i$ are the connected components.
$S_i$ are not necessarily homeomorphic to each other and the action of $\widetilde\tau_1$ on each $S_i$ is not necessarily all equal.
Each $S_i$ $(1\leq i\leq m)$ is induced from $d_{i}$-fold branched covering map, so it has one or two boundary components.
Then the embedding $\phi:B_n\hookrightarrow\Gamma_{g,k}$ is decomposed into
$$B_n\xrightarrow{\Delta}\oplus_{i=1}^{m} B_n\xrightarrow{\oplus\phi_i'}\oplus_{i=1}^{m}\Gamma_{g_i,k_i}^{(n)}\xrightarrow{\gamma}\Gamma_{g,k}.$$
Here $\Delta$ is the diagonal map, each $\phi_i'$ is the map obtained by $d_i$-fold covering (the same map as dealt in the connected case), $\gamma$ is the map recovering the surfaces $I_i$.
Note that $k=\sum_{i=1}^{m}k_i$ and each $k_i$ is the number of boundary components of the surface obtained by $d_i$-fold covering.

Let $d = \text{lcm}(d_{1}, \cdots ,d_{m})$.
For a regular embedding $\phi : B_{n} \rightarrow \Gamma_{g,k}$, by taking $\mathcal{C}_{i} = \conf_{di}(D)$ and $\mathcal{S}_{i} = \mathcal{M}_{g(i), k}$, we may give the operad action as before (as given in section 2) to get a $\mathcal{D}$-algebra map $\Phi : \mathcal{C} \rightarrow \mathcal{S}$. Thus, the induced map $\B\phi:\conf_{di}(D)\rightarrow\mathcal{M}_{g(i),k}$ preserves the actions of the little 2-disks operad. 
This proves the homology triviality.
\end{proof}

\section*{Acknowledgments}
The third author was supported by the Korean National Research Fund NRF-2020R1F1A1A01071639.

\printbibliography

\end{document}